\documentclass[a4paper,11pt]{article}
\usepackage[utf8]{inputenc}
\usepackage[T1]{fontenc}
\usepackage[english]{babel}
\usepackage{amsfonts,dsfont}
\usepackage{fullpage,units}
\usepackage{amsmath,amsthm,amsfonts,amssymb,upgreek,upgreek}
\usepackage{graphicx}

\usepackage{mathptmx}

\newtheorem{thm}{Theorem}[section]
\newtheorem{prop}[thm]{Proposition}
\newtheorem{lem}[thm]{Lemma}

\newtheorem{remarq}[thm]{Remark}

\newcommand{\ve}{\varepsilon}
\newcommand{\Hm}{X^{(\upgamma)}}

\newcommand{\Um}{X^{(\mathtt e)}}
\newcommand{\Xm}{X^{(\varphi)}}
\newcommand{\Wm}{W}
\newcommand{\Cdot}{{\mbox{\tiny\textbullet}}}

\begin{document}

% En-t{\^e}te
\title{\Large\bf{{Existence and asymptotic behaviour of some time-inhomogeneous diffusions} }}
\author{Mihai Gradinaru and  Yoann Offret\\\\
{\small Institut de Recherche Math{\'e}matique de Rennes, 
Universit{\'e} de Rennes 1},\\ 
{\small Campus de Beaulieu, 35042 Rennes Cedex, France}\\
{\tt{\small \{Mihai.Gradinaru,Yoann.Offret\}@univ-rennes1.fr}}}
\date{}
\maketitle

% R{\'e}sum{\'e}
{\small\noindent {\bf Abstract: } Let us consider a solution of a one-dimensional stochastic differential equation driven by a standard Brownian motion with time-inhomogeneous drift coefficient 
$\rho\,{\rm sgn}(x)|x|^\alpha/t^\beta$. This process can be viewed as a Brownian motion evolving in a potential, possibly singular, depending on time. We prove results on the existence and uniqueness of solution, study its asymptotic behaviour and made a precise description, in terms of parameters $\rho,\alpha$ and $\beta$, of the recurrence, transience and convergence. More precisely, asymptotic distributions, iterated logarithm type laws and rates of transience and explosion are proved for such processes.}\\

{\small\noindent {\bf R{\'e}sum{\'e} :} Nous consid{\'e}rons la solution d'une {\'e}quation diff{\'e}rentielle stochastique, dirig{\'e}e par un mouvement brownien lin{\'e}aire standard, dont le terme de d{\'e}rive varie avec le temps 
$\rho\,{\rm sgn}(x)|x|^\alpha/t^\beta$. Ce processus peut {\^e}tre vu comme un mouvement brownien {\'e}voluant dans un potentiel d{\'e}pendant du temps,  {\'e}ventuellement singulier. Nous montrons des r{\'e}sultats
d'existence et d'unicit{\'e} et nous {\'e}tudions le comportement asymptotique de la solution. Les propri{\'e}t{\'e}s de r{\'e}currence ou de transience de cette diffusion sont d{\'e}crites en fonction des param{\`e}tres  $\rho,\alpha$ et $\beta$, 
et nous donnons les vitesses de transience et d'explosion. Des r{\'e}sultats de convergence en loi et des lois de type logarithme it{\'e}r{\'e} sont {\'e}galement obtenus.}\\

% Mots clefs

% Classification AMS

%%%%%%%%%%%%%%%%%%%%%%%%%%%%%%%%%%%%%%%%%%%%%%%%%%%%%%%%%%%%%%%%%%%%%%%%%%%%%%%%%%%%%%%%%%%%%%%%%%%%%%%%%%%%%%%%%%%%%%%%%%%%%%%%%%%%%%%%%%%%%%%%%%%%%%%%%%%%%%%%%%%%%%%%%%%%%%%%%%%%%%%%%%%%%%%%%%%%%%%%%%%%%%%%%%%%%%%%%%%%%%%%%%%%%%%%%%%%%%%%%%%%%%%%%%%%%%%%%

\section{Introduction}

Let $X$ be a one-dimensional process describing a Brownian motion dynamics in a moving, possibly singular, potential 
$V_{\rho,\alpha,\beta}\,$:
\begin{equation}\label{process}
d X_t=d B_t-\frac{1}{2}\partial_{x}V_{\rho,\alpha,\beta}(t,X_t)\,d t,\quad X_{t_{0}}=x_{0},
\end{equation}
with,
\begin{equation}\label{potential}
V_{\rho,\alpha,\beta}(t,x):=\frac{-2\rho}{\alpha+1}\frac{|x|^{\alpha+1}}{t^{\beta}},\;\;\mbox{if}\;\; \alpha\neq -1\quad\mbox{and}\quad V_{\rho,\alpha,\beta}(t,x):=\frac{-2\rho\,\log{|x|}}{t^{\beta}},\;\;\mbox{if}\;\; \alpha=-1,
\end{equation}
where $B$ denotes a standard linear Brownian motion, $t_{0}>0$ and $x_{0},\rho,\alpha,\beta$ are some real constants. 
In this paper we shall study the asymptotic behaviour of such process. More precisely, our main goal is to give conditions which characterise the recurrence, transience and convergence in terms of parameters $\rho,\alpha$ and $\beta$. Here are the natural  questions one can ask: does there exist pathwise unique strong solution $X$ for equation (\ref{process})?
is this solution $X$ recurrent or transient? does there exist a well chosen normalisation of $X$ to ensure that the normalised process converges in distribution or almost surely? is it possible to obtain pathwise largest deviations of $X$, for instance iterated logarithm type law?

Questions as the last two ones are treated in \cite{AM,AW,GS} for different equations having some common features with (\ref{process}). For instance, Gihman and Skorohod in \cite{GS}, Chap. 4, \S17, consider the following equation
\begin{equation}\label{GS-Appleby}
dY_{t}=dB_{t}+d(Y_t)dt,\,\mbox{ with }\, d(y)\underset{|y|\to\infty}{\sim}\rho\,|y|^\alpha,\quad \rho>0\,\mbox{ and}\,-1<\alpha<1.
\end{equation}
Under additional assumptions, one proves that $Y_t$ is transient and asymptotically behaves as a solution of the deterministic underlying dynamical system, that is 
\begin{equation*}
Y_t\underset{t\to\infty}{\sim} h_t\quad\mbox{a.s.,}\quad\mbox{with}\quad h^\prime_t=d(h_t).
\end{equation*}
Equation (\ref{GS-Appleby}) is also considered by Appleby and Wu \cite{AW} with particular $\alpha=-1$. 
Its study is related to the Bessel process and the situation is more difficult. One proves that $Y_t$ 
satisfies the iterated logarithm law and recurrence or transience depends on the position of $\rho$ with respect 
to $1/2$. Appleby and Mackey \cite{AM} study the following damped stochastic differential equation
\begin{equation}\label{Appleby}
dY_{t}=\sigma(t)dB_{t}+d(Y_t)dt,\,\mbox{ with }\,d(y)\underset{y\to 0}{\sim}\rho\,{ \rm sgn}(y)|y|^\alpha,\,\rho<0\,\mbox{ and }\,\alpha>1.
\end{equation}
Here the diffusion coefficient $\sigma\in{\rm L}^2$ is such that $\sigma(t)\downarrow0$, as $t\to\infty$.
It is proved that $Y_t$ converges almost surely to 0 with polynomial rate. 
We will see that equation (\ref{Appleby}) is connected to equation (\ref{process}) by performing a suitable change of time.

For time-homogeneous stochastic differential equations driven by a one-dimensional Brownian motion, there exist precise criteria for recurrence or transience (see, for instance, Kallenberg \cite{Ka}, Chap. 23), or explosion (see, for instance, Ikeda and Watanabe \cite{IW}, Chap. VI, \S3), using the scale function. Some of these criteria are extended to the time-inhomogeneous situation for dimension greater or equal than two in Bhattacharya and Ramasubramanian \cite{BR}. 
Unfortunately, the results in \cite{BR} do not apply to equation (\ref{process}), even it is stated that the method can be adapted to the one-dimensional case. Recall also that there exist some general results on recurrence or transience (see, for instance, Has'minskii \cite{Ha}, Chap. III), and explosion (see for instance Narita \cite{Na} or Stroock and Varadhan \cite{SV}, Chap. 10), based on the construction of some convenient Lyapunov functions. 
However, for equation (\ref{process}), the construction of such functions seems to be more delicate. 

Equation (\ref{process}) can be also viewed as a continuous counterpart of a discrete time model considered recently by Menshikov and Volkov \cite{MV}.
Indeed, the discrete time process studied in \cite{MV} is a random walk on the real positive half line such that 
\begin{equation*}
\mathds{E}(X_{t+1}-X_{t}\mid X_{t}=x)\underset{t\to\infty}{\sim}\rho\frac{x^{\alpha}}{t^{\beta}}.
\end{equation*} 
The authors establish when the  process is recurrent or transient for certain values of parameters $\rho,\alpha, \beta$,  give the answer to a open question concerning the Friedman's urn model (see Freedman \cite{Fre}), and present some open problems. Their approach is based on a precise study of some 
submartingales and supermartingales.

Contrary to the discrete time model one firstly needs to study the existence, uniqueness and explosion of solutions for (\ref{process}). In the present paper different situations are distinguished, following the values of $\rho$ and $\alpha$, and existence and uniqueness are proved. We point out 
that, when $\alpha<0$, the existence of a solution is not obvious, since the drift has a singularity. 
For the time-homogeneous case, a solution to this problem is given by Cherny and Engelbert in \cite{CE}, by using the scale function. 
These ideas do not apply to one-dimensional time-inhomogeneous stochastic differential equations, and this is 
the main difficulty of this part of our paper. Our idea is to use an appropriate change of time, taking full 
advantage of the scaling property of the Brownian motion, of the Girsanov transformation, but also of the classification of isolated singular points in \cite{CE}. These different tools, adapted to continuous time models, also allow to answer the question of explosion of the solution when $\alpha>1$. As an example, we point out that, when $2\beta>\alpha+1$, 
the solution explodes in finite time with a positive probability, but not almost surely.

\begin{figure}[!ht]
 \centering
 \includegraphics[width=13cm,height=3.3cm]
{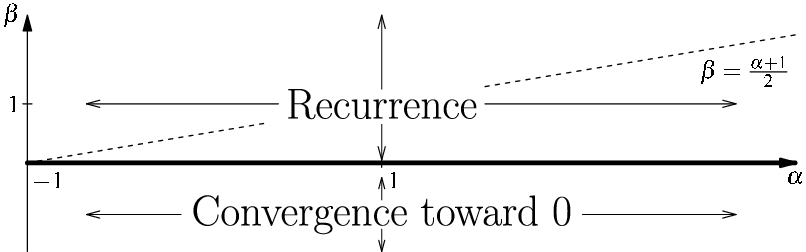}
%{phastransattract-1bis.ps}
 \caption{Phase transition in the attractive case $\rho<0$}
 \label{figattract}
 \end{figure}

Another goal of the present paper is to describe, for all values of parameters $\rho,\alpha,\beta$, the recurrent or the transient feature of the solution, but also its convergence. 
We present in Figure \ref{figattract} the diagram of phase transition that we obtain in the attractive case $\rho<0$.
\noindent
Note that $\alpha\leq-1$ is not allowed, since in that case, any solution is only defined up to the time of first reaching 0 (which is finite almost surely), and cannot be continued after it has reached this point. The critical line separating the two phases ({\sl recurrence} and {\sl convergence toward 0}) is $\beta=0$ and on this line the process is recurrent. The line $2\beta=\alpha+1$ could be called subcritical, in the sense that, the rate of the asymptotic behaviour is different on both sides. As for the proof of the existence, we use a suitable scaling transformation to obtain the asymptotic distribution of $X$ and its pathwise largest deviations, under a convenient normalisation. In fact, we show that the asymptotic behaviour of the process is strongly connected to the paths, and to the stationary distribution, of an ergodic diffusion. For example, when $2\beta<\alpha+1$, if $\varphi$ is the positive solution of $\varphi^\prime(t)=\varphi(t)^{\frac{2\beta}{\alpha+1}}$
 , then
\begin{equation*}
 \frac{X_{\varphi(t)}}{\sqrt{\varphi^\prime(t)}}\quad\quad \mbox{"behaves as"}\quad\quad H_t=B_t+\int_{0}^t\rho\,{\rm sgn}(H_s)|H_s|^\alpha ds. 
\end{equation*}
We obtain the convergence in distribution of $X_t / t^{\frac{\beta}{\alpha+1}}$ to the stationary distribution of $H$, and also its pathwise largest deviation. In particular, when $\beta<0$, we get the so-called polynomial stability of $X$. Furthermore, note that, if we set $Y_t:=X{(\phi_t)}$, with $\phi_t:=t^{\frac{1}{1-\beta}}$, then $Y_t$ satisfies the damped stochastic differential equation (\ref{Appleby}). We prove similar results as in  
\cite{AM} under slightly different hypothesis, and we obtain sharp rates of convergence.

\begin{figure}[!ht]
 \centering
 \includegraphics[width=13cm,height=3.3cm]
{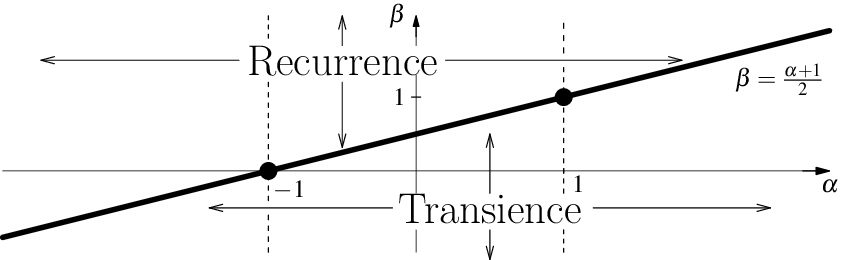}
%{phastransrepuls-1.ps}
 \caption{Phase transition in the repulsive case $\rho>0$}
 \label{figrepuls}
 \end{figure}

We present in Figure \ref{figrepuls} the diagram of phase transition that we obtain in the repulsive case $\rho>0$.
When $\alpha>1$ and $2\beta\leq\alpha+1$ the explosion time is almost surely finite. The critical curve is composed 
from two half-lines, $\beta=0$, when $\alpha\leq -1$, and $2\beta=\alpha+1$, when $\alpha\geq-1$,
and the process is either recurrent or transient. We prove similar results as in \cite{MV}, and again we obtain sharp rates of convergence. On the critical curve one needs to distinguish two particular points $(-1,0)$ and $(1,1)$, because these are the only cases where recurrence and transience depend on the position of $\rho$ with respect to $1/2$. 
$(\alpha,\beta)=(-1,0)$ corresponds to the well known Bessel process, whereas $(\alpha,\beta)=(1,1)$ is a continuous time counterpart of the Friedman's urn model. In the latter case, we obtain similar results as in \cite{Fre} and \cite{MV}, concerning recurrence and transience, but also regarding the asymptotic distribution and the pathwise largest deviations. For the other points of the critical curve, the process is recurrent. We point out that this is an open problem in \cite{MV}. 
The lines $\alpha=-1$ and $\alpha=1$ could be called subcritical, in the sense that, the behaviour of the process is slightly different to the right or left. In particular, the domain of recurrence depends on $\alpha$. The proof of recurrence is based on the same ideas as for the attractive case: by an appropriate scaling transformation of $X$ 
we associate an ergodic diffusion, whose asymptotic behaviour is easier to obtain. 
For instance, when $2\beta>\alpha+1$ and $-1<\alpha<1$, we show that
\begin{equation*}
 \frac{X_{e^{t}}}{e^{\frac{t}{2}}}\quad\quad \mbox{"behaves as"}\quad \quad U_t=B_t-\int_0^t \frac{U_s}{2} ds.
\end{equation*}
We get that $X$ behaves as a standard Brownian motion: it satisfies the iterated logarithm law and $X_t/\sqrt{t}$ converges in distribution to a standard Gaussian random variable. Roughly speaking, this means that 
the drift is asymptotically negligible compared to the noise. Concerning the proof of the transient case, when $\alpha<1$, the tools are similar to those used in \cite{MV}. We obtain similar results as in \cite{GS}, for equation (\ref{GS-Appleby}), and we show that $X$ behaves as a solution of the deterministic underlying dynamical system, that is
\begin{equation*}
|X_t|\underset{t\to\infty}{\sim} |h_t|\quad\mbox{a.s.,}\quad\mbox{with}\quad h_t^\prime=\rho\,{\rm sgn}(h_t)\frac{{|h_t|^\alpha}}{{t^\beta}}.
\end{equation*}

Some results in the present paper could be obtained, with similar arguments, for a general potential $V$, 
under convenient assumptions, for instance, when $\partial_xV(t,x)=-2 f(t)g(x)$ with $f(t)\sim_{t\to \infty} t^{-\beta}$ and $|g(x)|\sim_{|x|\to\infty}\rho|x|^{\alpha}$. These results will be presented elsewhere.
The case of a multiplicative noise seems more difficult.
Another interesting situation is obtained when one replaces the Brownian motion by a (stable) L{\'e}vy process, 
and it is object of some works in progress. Some methods in the present paper can be used in the study of 
a time-inhomogeneous diffusion in random environment of the form $V(t,x)=t^{-\beta}\,W(x)$, with $W$ a self-similar process  (for instance, a Brownian motion). This situation it is also object of some works in progress.

The paper is organised as follows: in the next section we introduce the scaling transformations and list the 
associated equations associated to some particular transformations. In Section 3 we perform the complete study 
of the existence, uniqueness and explosion of the solutions for equation (\ref{process}). Section 4 is devoted to 
a systematic study of the asymptotic behaviour of the solutions. Three cases are considered: on the critical line 
$2\beta=\alpha+1$, above and under this line. Proofs of some technical results are given in the Appendix. 

%%%%%%%%%%%%%%%%%%%%%%%%%%%%%%%%%%%%%%%%%%%%%%%%%%%%%%%%%%%%%%%%%%%%%%%%%%%%%%%%%%%%%%%%%%%%%%%%%%%%%%%%%%%%%%%%%%%%%%%%%%%%%%%%%%%%%%%%%%%%%%%%%%%%%%%%%%%%%%%%%%%%%%%%%%%%%%%%%%%%%%%%%%%%%%%%%%%%%%%%%%%%%%%%%%%%%%%%%%%%%%%%%%%%%%%%%%%%%%%%%%%%%%%%%%%%%%%%%%%%%%%%%%%%%%%%%%%%%%%%%%%%%%%%%%%%%%%%%%%%%%%%%%%%%%%%%%%%%%%%

\section{Scaling transformation and associated equations}
\setcounter{equation}{0}

We shall study equation (\ref{process}) in its equivalent form:
\begin{equation}\label{eqsign}
dX_t=dB_t+\rho\,{\rm sgn}(X_t)\frac{|X_t|^{\alpha}}{t^\beta}dt,\quad 
X_{t_{0}}=x_{0},
\end{equation}
$B$ being a standard Brownian motion defined on a filtered probability space $(\Omega,{\mathcal F},({\mathcal F}_{t})_{t\geq 0},\mathds{P})$. 
By symmetry of the equation and by the usual scaling transformation, we can assume without loss of generality that  $x_0\geq 0$ and $t_0=1$. We will keep these assumptions all along the paper.

We begin by defining a transformation of equation (\ref{eqsign}) which takes full advantage of the scaling property of the Brownian motion $B$ and the homogeneous properties of the drift $d(t,x):=\rho\,{\rm sgn}(x)|x|^\alpha/t^\beta$, that is for any $\lambda,\mu>0$, 
\begin{equation*}
(t\mapsto B_{\lambda t})\overset{\mathcal L}{=}(t\mapsto \lambda^{\frac{1}{2}}B_{t})\quad\mbox{and}\quad d(\mu t,\lambda x)= {\lambda^\alpha}\mu^{-\beta} d(t,x).
\end{equation*}
This transformation will provide some important equations related to our problem and it will be useful later to study the existence, the uniqueness and the asymptotic behaviour of solutions of equation (\ref{eqsign}).

%%%%%%%%%%%%%%%%%%%%%%%%%%%%%%%%%%%%%%%%%%%%%%%%%%%%%%%%%%%%%%%%%

\subsection{Scaling transformation}\label{scalingsection}

For any $T\in(0,\infty]$, let ${\rm \overline C}([0,T))$ be the set of functions $\omega : [0,T)\rightarrow
\mathds R \cup {\{\Delta\}} $ 
such that there exists a time $\tau_e(\omega)\in(0,T]$ (called the killing time of $\omega$) such that $\omega$ is continuous on $[0,\tau_e(\omega))$ and $\omega=\infty$ on $[\tau_e(\omega),T)$. 
We set $\Omega:={\rm \overline C} ([1,\infty))$ and $\Omega^\ast:={\rm \overline C}([0,t_{1}))$, with $t_1\in(0,\infty]$. For every ${\rm C}^{2}$-diffeomorphism (change of time)
$\varphi:[0,t_{1})\to[1,\infty)$  we introduce the scaling transformation $\Phi_\varphi:\Omega\to\Omega^\ast$ given by
\begin{equation}\label{scalingt}
\Phi_\varphi(\omega)(s):=\frac{\omega(\varphi(s))}{\sqrt{\varphi^\prime(s)}},\quad\mbox{with}\quad s\in[0,t_{1}),\,\omega\in \Omega.
\end{equation}
\begin{prop}\label{changetime}
The scaling transformation $\Phi_\varphi$ induces a bijection between weak solutions (possibly explosive) of equation (\ref{eqsign}) and weak solutions (possibly explosive)  of equation
\begin{equation}\label{eqtimechange}
d\Xm_s=d\Wm_s + \rho\,\frac{\varphi^\prime(s)^{\frac{\alpha+1}{2}}} 
{\varphi(s)^\beta}\,{\rm sgn}(\Xm_s)|\Xm_s|^{\alpha}\,ds - \frac{\varphi^{\prime\prime}(s)}
{\varphi^\prime(s)}\,\frac{\Xm_s}{2}\,ds,\quad \Xm_0=\frac{x_0}{\sqrt{\varphi^\prime(0)}}.
\end{equation}
Here $\{\Wm_s:s\in[0,t_{1})\}$ denotes a standard Brownian motion.
More precisely,
\begin{itemize}
\item[i)] if $(X,B)$ is a solution of equation (\ref{eqsign}) then $(\Xm,\Wm)$ is a solution of equation (\ref{eqtimechange}) where
\begin{equation}\label{PL1}
\Xm=\Phi_\varphi(X)\quad\mbox{and}\quad \Wm_t:=\int_0^t\frac{dB{(\varphi(s))}}{\sqrt{\varphi^\prime(s)}};
\end{equation}
\item[ii)] if $(\Xm,\Wm)$ is a solution of equation (\ref{eqtimechange})  then $(X,B)$ is a solution of equation (\ref{eqsign}) where
\begin{equation}\label{PL2}
X=\Phi_\varphi^{-1}(\Xm)\quad\mbox{and}\quad B_t-B_1:=\int_{1}^t\sqrt{(\varphi^\prime\circ\varphi^{-1})(s)}\,d\Wm(\varphi^{-1}(s)).
\end{equation} 
\end{itemize}
It follows that uniqueness in law, pathwise uniqueness or strong existence holds for equation (\ref{eqsign}) if, and only if, it holds for equation (\ref{eqtimechange}).
\end{prop}

\noindent
{\bf Proof.} 
Let $(X,B)$ be a solution of equation (\ref{eqsign}). By using P. L{\'e}vy's characterisation theorem of Brownian motion, we can see that $\Wm$ defined in (\ref{PL1}) is a standard Brownian motion. Moreover, 
by performing the change of variable $t:=\varphi(s)$ in (\ref{eqsign}), we get
\begin{equation*}
X_{\varphi(s)}-X_{\varphi(0)}=\int_{0}^{s}\sqrt{\varphi^\prime(u)}\,d\Wm_{u}+\rho\int_{0}^{s}{\rm sgn}\left(X_{\varphi(u)}\right) \frac{|X_{\varphi(u)}|^\alpha}{\varphi(u)^\beta}\,\varphi^\prime(u)\,du.
\end{equation*}
By the integration by parts formula written in its differential form we obtain 
\begin{equation*}
d\left(\frac{X_{\varphi(s)}}{\sqrt{\varphi^\prime(s)}}\right)
=d\Wm_s+\rho\,\frac{\varphi^\prime(s)^{\frac{\alpha+1}{2}}} 
{\varphi(s)^\beta}\,{\rm sgn}\left(\frac{X_{\varphi(s)}}{\sqrt{\varphi^\prime(s)}}\right)
\left|\frac{X_{\varphi(s)}}{\sqrt{\varphi^\prime(s)}}\right|^{\alpha}ds
-\frac{X_{\varphi(s)}}{\sqrt{\varphi^\prime(s)}}\,
\frac{\varphi^{\prime\prime}(s)}{2\varphi^\prime(s)}\,ds.
\end{equation*}
We conclude that equation (\ref{eqtimechange}) is satisfied by $(\Xm,W)$. The proof of (\ref{PL2}) is similar by noting that 
$\Phi_\varphi$ is a bijection and its inverse function is given by 
\begin{equation*}
\Phi_{\varphi}^{-1}(\omega)(s)=\sqrt{\varphi^\prime\circ\varphi^{-1}(s)}\,\omega(\varphi^{-1}(s)),\quad\mbox{with}\quad  s\in[1,\infty),\,\omega\in\Omega^*.
\end{equation*}
The final remark is a simple application of parts i) and ii).\hfill$\Box$

%%%%%%%%%%%%%%%%%%%%%%%%%%%%%%%%%%%%%%%%%%%%%%%%%%%%%%%%%%%%%%%%%%

\subsection{Two particular transformations}\label{scaling}

We give here two scaling transformations which produce at least one time-homogeneous coefficient among
the two terms of the drift in (\ref{eqtimechange}). We also introduce some equations related 
to (\ref{eqsign}) which will be useful later in our study. For simplicity, we will keep the notation 
$W$ for a standard Brownian motion which can be different from the process employed in Proposition 
\ref{changetime}. 

\subsubsection{Exponential scaling transformation}\label{exponentialscaling}

The transformation (\ref{scalingt}) associated to the exponential change of time $\varphi_{\mathtt e}(t):=e^t$, and denoted by $\Phi_{\mathtt e}$, is given by
\begin{equation*}
\Phi_{\mathtt e}(\omega)(s)=\frac{\omega(e^s)}{e^{\nicefrac{s}{2}}},\quad\mbox{with}\quad s\in[0,\infty),\,\omega\in\Omega.
\end{equation*}
The process $X^{(\mathtt e)}:=\Phi_{\mathtt e}(X)$ satisfies equation (\ref{eqtimechange}) which can be written
\begin{equation}\label{eqtimechangeexp}
d\Um_s=dW_s-\frac{\Um_s}{2}ds + \rho\,e^{\left(\frac{\alpha+1}{2}-\beta\right)s}\,{\rm sgn}(\Um_s)|\Um_s|^{\alpha}ds,
\quad \Um_0=x_0.
\end{equation}
On the one hand, if we leave out the third term in the right hand side of equation (\ref{eqtimechangeexp}), we obtain the classical equation of the Ornstein-Uhlenbeck process: 
%, that is
\begin{equation}\label{OU}
dU_s=dW_s -\frac{U_s}{2}ds,\quad U_0=x_0.
\end{equation}
Note that equation (\ref{OU}) is a particular case of equation (\ref{eqsign}) with parameters 
$\rho=-1/2$, $\alpha=1$ and $\beta=0$. On the other hand, when $\alpha=-1$, by  Ito's formula we can see that $Y:=X^2$ satisfies 
\begin{equation}\label{besselcarre}
dY_t=2\sqrt{Y_t}\,dW_t+\Big(\frac{2\rho}{t^\beta} +1\Big)\,dt,\quad Y_0=x_0^2,\quad Y\geq 0. 
\end{equation}
This process can be viewed as a square Bessel process whose dimension depends on time. Clearly, when
$\beta=0$, this process is the classical square Bessel process $R$ of dimension $2\rho+1$ and it satisfies
\begin{equation}\label{squarebessel}
dR_t=2\sqrt{R_t}\,dW_t+(2\rho +1)\,dt,\quad R_0=x_0^2,\quad R\geq 0.
\end{equation} 
Furthermore, the process $R^{(\mathtt e)}:=\Phi_{\mathtt e}(R)$ satisfies
\begin{equation}\label{eqtimechangebesselcarre}
dR^{(\mathtt e)}_t = 2\sqrt{R^{(\mathtt e)}_t}\,dW_t + \Big( 2\rho  +1 -\frac{R^{(\mathtt e)}_t}{2}\Big)\,dt,\quad R^{(\mathtt e)}_0 = x_0^2,\quad R^{(\mathtt e)}\geq 0.
\end{equation}

\subsubsection{Power scaling transformation}

Assume that $\alpha\neq -1$ and consider the Cauchy problem:
\begin{equation}\label{gamma}
\varphi_\upgamma^\prime(s)=\varphi_\upgamma(s)^\gamma,\quad\varphi_\upgamma(0)=1,\quad\mbox{ with }\quad
\gamma:=\frac{2\beta}{\alpha+1}.
\end{equation}
There exists a unique maximal solution $\varphi_\upgamma\in{\rm C}^2([0,t_{1});[1,\infty))$ and we can see that
\begin{equation*}
\varphi_\upgamma(s)=(1+(1-\gamma)s)^{\frac{1}{1-\gamma}},\quad\mbox{when}\quad\gamma\neq 1 \quad\mbox{and}\quad \varphi_\upgamma(s)=e^{s},\quad\mbox{when}\quad \gamma=1,
\end{equation*}
with $t_1=\infty$, when $\gamma\in(-\infty,1]$, and $t_1={1}/{(\gamma-1)}$, when $\gamma\in(1,\infty)$. 
The transformation (\ref{scalingt}) associated to this change of time will be denoted $\Phi_\upgamma$, and is  given by
\begin{equation*}
\Phi_\upgamma(\omega)(s)=\frac{\omega(\varphi_\upgamma(s))}{\varphi_\upgamma(s)^{\frac{\gamma}{2}}},\quad\mbox{with}\quad s\in[0,t_1),\,\omega\in\Omega.
\end{equation*}
The process $\Hm:=\Phi_{\upgamma}(X)$ satisfies equation (\ref{eqtimechange}) which can be written 
\begin{equation}\label{eqtimechangegamma}
d\Hm_s=dW_s+\rho\,{\rm sgn}(\Hm_s)|\Hm_s|^{\alpha}\,ds - {\gamma}\varphi_\upgamma^{\gamma-1}(s) \frac{\Hm_s}{2}\,ds,
\quad \Hm_0=x_0,\quad s\in[0,t_1).
\end{equation}
i) If $\gamma\in(-\infty,1)$, equation (\ref{eqtimechangegamma}) takes the form:
\begin{equation}\label{eqtimechangegamma3}
d\Hm_s=dW_s + \rho\,{\rm sgn}(\Hm_s)|\Hm_s|^{\alpha}\,ds -\frac{\gamma\,\Hm_s}{2(1+(1-\gamma)s)}ds,
\quad\Hm_0=x_0,\quad s\in[0,\infty).
\end{equation}
ii) If $\gamma\in(1,\infty)$, equation (\ref{eqtimechangegamma}) takes the form:
\begin{equation}\label{eqtimechangegamma2}
d\Hm_s=dW_s + \rho\,{\rm sgn}(\Hm_s)|\Hm_s|^{\alpha}\,ds -\delta \,\frac{\Hm_s}{t_{1}-s}\,ds ,
\quad\Hm_0=x_0,\quad s\in[0,t_1),
\end{equation}
with
\begin{equation*}%\label{delta_t1}
t_1=\frac{1}{\gamma-1}\quad\mbox{and}\quad \delta:=\frac{\gamma}{2(\gamma-1)}.
\end{equation*}
iii) If $\gamma=1$, equation  (\ref{eqtimechangegamma}) takes the form:
\begin{equation}\label{criticalline}
dZ_s=dW_s + \Big(\rho\,{\rm sgn}(Z_s)|Z_s|^{\alpha} -\frac{Z_s}{2}\Big)ds,
\quad Z_0=x_0,\quad s\in[0,\infty).
\end{equation}
Note that the transformations $\Phi_{\mathtt e}$ and $\Phi_\upgamma$ coincide when $\gamma=1$. Finally, let us introduce another two stochastic differential equations related to 
(\ref{eqtimechangegamma}). First, we leave out the third term on the right hand side of (\ref{eqtimechangegamma}) 
and we get
\begin{equation}\label{undercriticallinehomogeneous}
dH_s=dW_s+\rho\,{\rm sgn}(H_s)|H_s|^{\alpha}ds,\quad H_0=x_0,\quad s\in[0,t_1).
\end{equation}
Note that the latter equation is nothing but (\ref{eqsign}) with $\beta=0$.
Second, we leave out the second term in the right hand side of (\ref{eqtimechangegamma2}) and we obtain 
\begin{equation}\label{deltapontmansuy}
db_s=dW_s -\delta\,\frac{ b_s}{t_1-s}\,ds,\quad b_0=x_0,\quad s\in[0,t_1).
\end{equation}
The process $b$ is the so-called $\delta$-Brownian bridge (see also \cite{M}) and it is the classical 
Brownian bridge when $\delta=1$. 

%%%%%%%%%%%%%%%%%%%%%%%%%%%%%%%%%%%%%%%%%%%%%%%%%%%%%%%%%%%%%%%%%%%%%%%%%%%%%%%%%%%%%%%%%%%%%%%%%%%%%%%%%%%%%%%%%%%%%%%%%%%%%%%%%%%%%%%%%%%%%%%%%%%%%%%%%%%%%%%%%%%%%%%%%%%%%%%%%%%%%%%%%%%%%%%%%%%%%%%%%%%%%%%%%%%%%%%%%%%%%%%%%%%%%%%%%%%%%%%%%%%%%%%%%%%%%%%%

\section{Preliminary study of solutions}
\setcounter{equation}{0}

In this section we study existence, uniqueness and explosion of solutions for equation (\ref{eqsign}).
For parameters $(\rho,\alpha,\beta)\in\big(\mathds R\times(-1,\infty)\times \mathds R\big)\cup\big([0,\infty)\times(-\infty,-1] \times \mathds R\big)$ and $x_0\in[0,\infty)$ we prove the existence of a time-inhomogeneous diffusion $X$ solution of equation (\ref{eqsign}), defined up to the explosion time, and taking its values in  $\mathds R$, provided $\alpha\in(-1,\infty)$, in $(0,\infty)$, provided $\alpha\in(-\infty,-1)$ 
and in $[0,\infty)$, provided $\alpha=-1$. We show that this diffusion can explode in finite time with positive probability when $(\rho,\alpha,\beta)
\in(0,\infty)\times(1,\infty)\times\mathds R$.

%%%%%%%%%%%%%%%%%%%%%%%%%%%%%%%%%%%%%%%%%%%%%%%%%%%%%%%%%%%%%%%%%

\subsection{Existence and uniqueness}

Existence and uniqueness for equation (\ref{eqsign}) are not obvious since  the drift 
could be singular in 0 and/or not time-homogeneous. However, with the help of transformation (\ref{scalingt}), the Girsanov transformation and the results on power equations in \cite{CE}, chap. 5, %(see Lemma \ref{girsanov} below) 
we reduce to the study of equation (\ref{undercriticallinehomogeneous}) when $\alpha\neq -1$. The following remark is stated only for later reference. 

\begin{remarq}\label{girsanov}
Assume that $\alpha\in \mathds{R}\setminus\{-1\}$. The Girsanov transformation induces a linear bijection between weak solutions (respectively nonnegative solutions or positive solutions) defined up to the explosion time of equation (\ref{eqtimechangegamma}) and weak solutions (respectively nonnegative solutions or positive solutions) defined up to the explosion time of equation (\ref{undercriticallinehomogeneous}).
\end{remarq}

%%%%%%%%%%%%%%%%%%%%%%%%%%%%%%%%%%%%%%%%%%%%%%%%%%%

\subsubsection{Locally integrable singularity : $\alpha>-1$}

In this case  $x\mapsto |x|^{\alpha}$ is locally integrable. As for equation (\ref{undercriticallinehomogeneous}), we show that there exists a pathwise unique strong solution $X$ to equation (\ref{eqsign}) defined up to the explosion time.

\begin{prop}\label{existence2}
If $\alpha\in(-1,\infty)$, $\beta,\rho\in\mathds R$, there exists a pathwise unique strong solution $X$ to equation (\ref{eqsign}) defined up to the explosion time.
\end{prop}

\noindent
{\bf Proof.} By using Proposition 2.2 in \cite{CE}, p. 28, there exists a unique weak solution $H$ to the time-homogeneous equation (\ref{undercriticallinehomogeneous}) defined up to the explosion time. Therefore, by using Remark \ref{girsanov}, there exists a unique weak solution $X^{(\upgamma)}$ to equation (\ref{eqtimechangegamma}) and, by using Proposition \ref{changetime}, there exists a unique weak solution $X$ to equation (\ref{eqsign}). Moreover, since pathwise uniqueness holds for equation (\ref{eqsign}) by using Proposition 3.2 and Corollary 3.4, Chap. IX in \cite{RY}, pp. 389-390, we get the conclusion.
Note that for the nonsingular case $\alpha\geq 0$, the coefficients of equation (\ref{eqsign}) are continuous, hence the present proposition can be obtained by usual techniques (localisation, Girsanov and Novikov theorems).\hfill$\Box$\\

Let us denote by $\mathcal L(\pm X)$ the distribution of the process $\pm X$. We shall say that a probability distribution $\mu$ is a {\sl mixture of distributions of $X$ and $-X$}, if there exists $\lambda\in[0,1]$ such that $\mu=\lambda \mathcal L(X)+(1-\lambda)\mathcal L(-X)$. Equivalently, there exists a discrete random variable $U\in\{-1,+1\}$, independent of $X$, such that $\mu=\mathcal L(UX)$.

\subsubsection{Nonlocally integrable singularity : $\alpha<-1$ and $\rho>0$}

Again, it suffices to study equation (\ref{undercriticallinehomogeneous}). We shall see that there exists a pathwise unique nonnegative solution $X$ to equation (\ref{eqsign}) and there are several strong Markov weak solutions when the process start at the singularity $x_0=0$. 

\begin{prop}\label{existence3} 
If $\alpha\in(-\infty,-1)$,
$\beta\in \mathds{R}$ and $\rho\in(0,\infty)$, there exists a pathwise unique nonnegative strong solution $X$ to equation (\ref{eqsign}).
Moreover,
\begin{enumerate}
\item[i)] if $x_0\in(0,\infty)$, $X$ is the pathwise unique strong solution and it is positive;
\item[ii)] if $x_0=0$, for all $t>1$, $X_t>0$ a.s. and the set of all weak solutions is the set of all distributions which are mixture of distributions of $X$ and $-X$.
\end{enumerate}
\end{prop}

\noindent
{\bf Proof.}
By using Theorem 3.5 in \cite{CE}, p. 66, there exists a unique nonnegative weak solution $H$ to equation (\ref{undercriticallinehomogeneous}) defined up to the explosion time. We deduce, following the same lines as in the proof of the previous proposition, that there exists a pathwise unique nonnegative strong solution to equation (\ref{eqsign}).
We point out that there is no uniqueness in law for equation (\ref{eqsign}) when $x_0=0$ and we cannot apply  directly Proposition 3.2, Chap. IX in \cite{RY}, p. 389, to prove the the pathwise uniqueness. However, we can check that a similar result to the cited proposition, whose the proof 
can be imitated, holds for nonnegative solutions.  

Moreover, if $x_0\in(0,\infty)$, any weak  solution of equation (\ref{undercriticallinehomogeneous}) is positive and we deduce that there exists a pathwise unique strong solution to equation (\ref{eqsign}) and this solution is positive. 

Finally, if $x_0=0$, the set of all weak  solutions of equation (\ref{undercriticallinehomogeneous}) is (by symmetry of the equation) the set of all distributions which are mixture of the distributions of $H$ and $-H$. We deduce the point ii) and the proof is done.\hfill$\Box$

\subsubsection{Bessel type case: $\alpha=-1$ and $\rho>0$}

In this case, Remark \ref{girsanov} does not hold and we perform a direct study of (\ref{eqsign}). We show that there exists a pathwise unique nonnegative strong solution $X$ to equation (\ref{eqsign}), which can be viewed as a Bessel process whose the dimension $2\rho t^{-\beta}+1$ depends on time. Note that it is possible that there exists different weak solutions (not necessarily Markovian). 

\begin{prop}\label{besselinhom}
If $\alpha=-1$, $\rho\in(0,\infty)$ and $\beta\in \mathds R$, there exists a pathwise unique nonnegative strong solution $X$ to equation (\ref{eqsign}). Moreover,
\begin{enumerate}
\item[i)] if $\rho\in\left[{1}/{2},\infty\right)$, $\beta\in(-\infty,0]$
and $x_0\in(0,\infty)$, $X$ is the pathwise unique strong solution and it is positive;
\item[ii)] if $\rho\in\left[{1}/{2},\infty\right)$, $\beta\in(-\infty,0]$
and $x_0=0$, the set of all weak solutions is the set of all distributions which are mixture of distributions of $X$ and $-X$ and $\forall\;t>0$, $X_t>0$ a.s.;
\item[iii)]  if $\beta\in(0,\infty)$  or if $(\rho,\beta)\in\left(0,{1}/{2}\right)\times (-\infty,0]$, we can construct different weak solutions to equation (\ref{eqsign}) and in the first case the set $\{t\geq 1 : X_t=0\}$ is unbounded a.s.
\end{enumerate}
\end{prop}

\noindent
{\bf Proof.} 
To begin with, it is not difficult to see that there exists a pathwise unique nonnegative strong solution $Y$ to equation (\ref{besselcarre}). This process can be viewed as the squared Bessel process having a time-dependent dimension $2\rho t^{-\beta}+1$. We shall prove that $X:=\sqrt{Y}$ is a nonnegative weak solution of equation (\ref{eqsign}). By applying Ito's formula, for all $t\geq 1$ and $\varepsilon>0$, 
\begin{equation}\label{bessitoepsilon}
(X_t^2+\varepsilon)^{\frac{1}{2}}=(x_0^2+\varepsilon)^{\frac{1}{2}}+\int_{1}^{t} 
\Big(\frac{X_s^2}{X_s^2+\varepsilon}\Big)^{\frac{1}{2}}dW_s
+\int_{1}^t\frac{\rho\,ds}{s^\beta (X_s^2+\varepsilon)^{\frac{1}{2}}} 
+ \int_{1}^t\frac{\ve\,ds}{2(X_s^2+\varepsilon)^{\frac{3}{2}}}.
\end{equation} 
We let $\ve \to 0$ in (\ref{bessitoepsilon}). Firstly, it is clear that 
\begin{equation*}
\lim_{\varepsilon\to 0}\int_{1}^t \Big(\frac{X_s^2}{X_s^2+\varepsilon}\Big)^{\frac{1}{2}}dW_s
=W_t-W_{1},\quad\mbox{ in probability.}
\end{equation*}
Secondly, by monotone convergence theorem, the third term in the right hand 
side of (\ref{bessitoepsilon}) converges a.s. We show that the limit is finite a.s. and that the fourth term converges toward $0$ in probability by comparison with a squared Bessel process. To this end, let us consider the pathwise unique nonnegative strong solution of
\begin{equation}\label{besselcompar}
Q_s= x_0+ W_s-W_1 + \int_{1}^s\frac{\rho_1}{Q_u}\,du,\quad s\geq 1,\quad\mbox{with}\quad \rho_1:=\inf\left\{\frac{\rho}{s^\beta}:s\in[1,t]\right\}>0.
\end{equation}
$Q$ is a classical Bessel process of dimension $2\rho_1+1$.
By using a comparison theorem (see Theorem 1.1, Chap. VI in [11], p. 437) and Ito's formula, we can see that for all $s\in[1,t]$, $X_s^2\geq Q_s^2$ and
\begin{equation*}
(Q_t^2+\varepsilon)^{\frac{1}{2}}=(x_0^2+\varepsilon)^{\frac{1}{2}}+\int_{1}^{t} 
\Big(\frac{Q_s^2}{Q_s^2+\varepsilon}\Big)^{\frac{1}{2}}dW_s
+\int_{1}^t\frac{\rho_1\, ds}{(Q_s^2+\varepsilon)^{\frac{1}{2}}} 
+  \int_{1}^t\frac{\ve\,ds}{2(Q_s^2+\varepsilon)^{\frac{3}{2}}}.
\end{equation*}
Since $Q$ is a solution of (\ref{besselcompar}) we obtain, by letting $\ve\to 0$ in the latter equality,
\begin{equation*}
\int_{1}^t\frac{\rho\,ds}{s^\beta X_s}\leq \int_{1}^t\frac{\rho_2}{Q_s}ds<\infty\quad\mbox{a.s.}\quad\mbox{with}\quad \rho_2:=\sup\Big\{\frac{\rho}{s^\beta}:s\in[1,t]\Big\}<\infty,
\end{equation*}
and also
\begin{equation*}
\lim_{\ve\to 0}\int_{1}^{t}\frac{\ve\,ds}{ (X_s^2+\varepsilon)^{\frac{3}{2}}}\leq 
\lim_{\ve\to 0}\int_{1}^{t}\frac{\ve\,ds}{ (Q_s^2+\varepsilon)^{\frac{3}{2}}}
=0,\quad\mbox{in probability.}
\end{equation*}
We get that $X$ is a nonnegative weak solution of (\ref{eqsign}).  Pathwise uniqueness is obtained by 
using the same arguments as in the proof of Propositions \ref{existence2} and \ref{existence3} and we deduce that there exists a pathwise unique nonnegative strong solution $X$ to equation (\ref{eqsign}). 
We proceed with the proof of points i)-iii) in the statement of the proposition.

Firstly, if $\rho\in[{1}/{2},\infty)$, $\beta\in(-\infty,0]$
and $x_0\in(0,\infty)$, the inequality $2\leq 2\rho t^{-\beta}+1$ holds for all $t\geq 1$ and we deduce that $X$ 
is positive by comparison with a Bessel process of dimension $2$.

Secondly, if $\rho\in[{1}/{2},\infty)$, $\beta\in(-\infty,0]$ and $x_0=0$, the same comparison can be used to see that every solution $ \tilde X$ of (\ref{eqsign}) satisfies $\tilde X^{2}_{t}\neq 0$, for all $t>1$ a.s. Let us introduce 
\begin{equation*}
\Omega^\pm:=\{\omega\in \Omega : \forall t>1,\,\pm\tilde X_t>0\}
\quad\mbox{and}\quad \mathds{P}^\pm:=\mathds{P}\left(\,\Cdot\,\big| \Omega^\pm\right).
\end{equation*}
For all $\ve>0$, $\Omega^\pm=\{\omega\in \Omega: \forall 1<t<1+\ve,\,\pm\tilde X_{t}>0\}\in\mathcal{F}_{t+\ve}$ and then $\Omega^\pm\in\mathcal{F}_{1+}$.
Therefore, the standard Brownian motion $\{B_t-B_{1}\}_{t\geq 1}$ under $\mathds P$ is again a standard Brownian motion under probabilities $\mathds{P}^\pm$. By uniqueness of the nonnegative weak solution and also, by symmetry, of the nonpositive solution of (\ref{eqsign}),  the distribution of $\tilde X$ under $\mathds{P}^\pm$ equals to the distribution of $\pm X$. The point ii) is then a simple consequence. 

Finally, if $\beta\in(0,\infty)$, for $t$ large enough we have $2\rho t^{-\beta}+1\leq \delta$, with $\delta\in(0,1)$.
By comparison with a Bessel process of dimension $\delta$, we get that the reaching time of $0$ 
is finite a.s. and the set $\{t>1 : X_t=0\}$ is unbounded a.s. Besides, if $X$ is a solution starting from $x_0=0$, $-X$ is also a solution. We deduce that different solutions could be constructed by gluing the paths of $X$ and $-X$ each time when the process returns in $0$. If $\rho\in(0,{1}/{2})$ and $\beta\in(-\infty,0)$, for all $s\in[1,(2\rho)^{1/\beta})$ and $t\in[1,s]$, $2\rho t^{-\beta}+1\leq 2\rho s^{-\beta}+1$. We deduce by comparison with a Bessel process of dimension $2\rho s^{-\beta}+1\in(1,2)$ that the reaching time of $0$ belongs to $\left[1,(2\rho)^{\nicefrac{1}{\beta}}\right)$ with a positive probability. Indeed, the reaching time of $0$ for a Bessel process of this dimension has a positive density with respect to the Lebesgue measure 
(with an explicit expression given, for instance, in \cite{GRVY}, p. 537). As in the preceding case, different solutions can be constructed.
\hfill$\Box$

\begin{remarq}\label{parameters}
By using similar methods as in Propositions \ref{existence2}, \ref{existence3} and \ref {besselinhom}, when $\alpha\leq -1$ and $\rho<0$, it can be proved that weak solutions of equation (\ref{eqsign}) are only defined up to the reaching time of $0$, which is finite a.s. and cannot be continued after this time. This case will be not considered since is out of range for the study of the asymptotic behaviour.
\end{remarq}

%%%%%%%%%%%%%%%%%%%%%%%%%%%%%%%%%%%%%%%%%%%%%%%%%%%%%%%%%%%%%%

\subsection{Explosion of solutions}

We show that $X$ explodes in finite time with positive probability if and only if $\alpha\in(1,\infty)$. More precisely, the explosion time $\tau_e$ of $X$ is finite a.s., provided $2\beta\leq \alpha+1$, and satisfies $\mathds P(\tau_e=\infty)\in(0,1)$, provided $2\beta>\alpha+1$.

\begin{prop}\label{timeexplosion} The explosion time $\tau_e$ of $X$ is infinite a.s. if $\rho \in(-\infty,0)$ or $\alpha\in(-\infty,1]$.\, It is  finite a.s. if $\rho\in(0,\infty)$, $\alpha\in(1,\infty)$ and $2\beta\in(-\infty, \alpha+1]$.
\end{prop}

\noindent
{\bf Proof.} Assume first that $\rho\in(-\infty,0)$ or $\alpha\in(-\infty,1 ]$. Let $F$ be a twice continuous differentiable nonnegative function such that $F(x):=1+x^2$ for all $|x|\geq 1$, $F(x)=1$ for all $x\in[1/2,1/2]$ and $F\geq 1$. For all $T\geq 1$, we denote $c_T$ the supremum of $LF$ on $[1,T]\times [-1,1]$, where  
$L$ is the infinitesimal generator of $X$ given by
\begin{equation}\label{ig}
L:= \frac{1}{2}\frac{\partial^2}{\partial x^2} + \rho\,{\rm sgn}(x)\frac{|x|^\alpha}{t^\beta}\frac{\partial}{\partial x} +\frac{\partial}{\partial t} .
\end{equation}
It is a simple calculation to see that for all $t\in[1,T]$ and $x\in\mathds R$,
\begin{equation*}
LF(t,x)\leq c_T+\lambda_T F(t,x)\leq (c_T+\lambda_T) F(t,x),\quad\mbox{with}\quad 
\lambda_T:= \sup_{1\leq t\leq T}(1+|\rho|t^{-\beta}).
\end{equation*}
By using Theorem 10.2.1 in \cite{SV}, p. 254, we deduce that the explosion time $\tau_e$ is finite a.s.

Finally, assume that $\rho\in(0,\infty)$, $\alpha\in(1,\infty)$ and $2\beta\in(-\infty, \alpha+1]$. By using Proposition 
\ref{changetime} it suffices to show that the solution $\Hm$ of equation (\ref{eqtimechangegamma}) explodes in finite time a.s.
Let us introduce $Q_s$ and $C_s$, the pathwise unique strong solutions of 
\begin{equation*}
dQ_s = 2\sqrt{Q_s}\,dW_s+\Big(
2\rho\,Q_s^{\frac{\alpha+1}{2}} - 
|\gamma|Q_s+1\Big)\,ds,
\quad Q_0=x_0^2,
\end{equation*}
and
\begin{equation*}
dC_s = 2\sqrt{C_s}\,dW_s+\Big(
2\rho\,C_s^{\frac{\alpha+1}{2}} - 
\gamma \varphi_{\upgamma}^{\gamma-1}(s)\,C_s+1\Big)\,ds,
\quad C_0=x_0^2.
\end{equation*}
By using Ito's formula, we can see that the square of $\Hm$ satisfies the latter equation and by weak uniqueness, we get that  
$C$ and $(\Hm)^2$ have the same distribution.
Moreover, since $\gamma={2\beta}/{(\alpha+1)}\leq 1$, we can see that $0\leq \varphi^{\gamma-1}_{\upgamma}\leq 1$. By comparison theorem, 
we get that  $0\leq Q_s\leq C_s$ a.s. Besides, by using Theorem 5.7 in \cite{CE}, p. 97, the explosion 
time of the time-homogeneous diffusion $Q$ is finite a.s. We deduce that the explosion time of $C$, and consequently that of $\Hm$, is finite a.s.\hfill$\Box$

\begin{prop}\label{timeexplosion2}
If $\rho\in(0,\infty)$, $\alpha\in(1,\infty)$ and $2\beta\in(\alpha+1,\infty)$,
\begin{equation}\label{explosiontime}
\mathds P(\tau_e=\infty)=\mathds{E}\Big[\exp\Big(
\int_0^{t_1} \rho\,{\rm sgn}(b_u)|b_u|^{\alpha}dW_u-\frac{1}{2}
\int_0^{t_1}\rho^2 |b_u|^{2\alpha}du
\Big)\Big]\in(0,1),
\end{equation} 
where $b$ denotes the weak solution of equation (\ref{deltapontmansuy}) and $\tau_e$ the explosion time of $X$. 
\end{prop}

\noindent
{\bf Proof.} Let $\Hm$ be the pathwise unique strong solution of equation (\ref{eqtimechangegamma2}) and $b$ be the pathwise unique strong solution of equation (\ref{deltapontmansuy}).
Recall that $\gamma={2\beta}/{(\alpha+1)}>1$ and $t_{1}={1}/{(\gamma-1)}$. 
Denote by $\eta_e$ the explosion time of $\Hm$ and note that a.s. $\eta_e  \in [0,t_1]\cup\{\infty\}$ and $\{\eta_e\geq t_1\}=\{\tau_e=\infty\}$. We need to show that $\mathds{P}(\eta_e\geq t_1)$ is equal to the right hand side of (\ref{explosiontime}) and belongs to $(0,1)$. First of all,  $b$ is a continuous process on $[0,t_1]$, with $b_{t_1}=0$ a.s.,  it is the so-called $\delta$-Brownian bridge (see Definition 1 in \cite{M}, p. 1022). By using the Girsanov transformation between $b$ and $X^{(\upgamma)}$, we can write for every integer $n\geq 1$, $s\in[0,t_{1}]$ and $A\in\mathcal{F}_s$,
\begin{equation*}
\mathds{E}\Big[\mathds{1}_A\big(\Hm_{\Cdot\wedge \eta_n}\big)\mathds{1}_{\{\eta_n> s\}}\Big]  
=\mathds{E}\Big[\mathds{1}_A\left(b_{\Cdot\wedge\sigma_n}\right)\mathcal{E}\left(s\wedge\sigma_{n}\right)\mathds{1}_{\{\sigma_n>s\}}
\Big],
\end{equation*}
where
\begin{equation*}
\eta_n:=\inf\{s\in[0,t_{1}) : |\Hm_{s}|\geq n\},\quad \sigma_n:=\inf\{s\in[0,t_{1}) : |b_s|\geq n\},
\end{equation*}
and
\begin{equation*}
\mathcal{E}(s):=\exp\Big(
\int_0^s \rho\,{\rm sgn}(b_u)|b_u|^{\alpha}dW_u-\frac{1}{2}
\int_0^s \rho^2 |b_u|^{2\alpha}du
\Big).
\end{equation*}
Letting $n\to\infty$, we obtain
\begin{equation}\label{girsa}
\mathds{E}\Big[\mathds{1}_A(\Hm)\mathds{1}_{\{\eta_e> s\}}\Big]= 
\mathds{E}\Big[\mathds{1}_A(b)\mathcal{E}(s)\Big].
\end{equation}
In particular, we have proved that for all $s\in[0,t_{1}]$, $\mathds{P}(\eta_e> s)  = 
\mathds{E}\left[\mathcal{E}(s)\right]$. Furthermore it is clear that $\mathds{P}(\tau_e=\infty)=\mathds{P}(\eta_e\geq t_1)\geq \mathds{E}\left[\mathcal{E}(t_1)\right]>0$. At this level we state a technical result which proof is 
postponed to the Appendix.
\begin{lem}\label{lemmepresquehorrible}
Assume that $\rho\in(0,\infty)$, $\alpha\in(1,\infty)$ and $2\beta\in(\alpha+1,\infty)$, and denote by $\eta_e\in[0,t_1]\cup\{\infty\}$ the explosion time of  $X^{(\upgamma)}$ (the weak solution of (\ref{eqtimechangegamma2})). Then $\mathds{P}(\eta_e=t_1)=0$.
\end{lem}

\noindent
We deduce from this lemma that $\mathds{P}(\tau_e=\infty)=\mathds P(\eta_e\geq t_1)=\mathds P(\eta_e>t_1)=\mathds E(\mathcal E(t_1))$ and the equality in (\ref{explosiontime}) is proved. It remains to show that $\mathds{P}(\tau_e=\infty)<1$. Recall that $\alpha\in(1,\infty)$ and let $a\in(1,\alpha)$. Set $g(x):=1\wedge|x|^{-a}$ and note that, for any $T>1$, we can choose $k\geq 1$, such that
$a(a-1)^{-1}=\int_0^\infty g(y)dy < k (T-1)$. 
Moreover, we can see that there exists a continuous differentiable odd function $f$, defined on $\mathds R$, 
vanishing only at $x=0$, such that $|f|\leq g$, and
\begin{equation*}
f(x):=k x,\quad x\in\left[-1/2k,1/2k\right],
\quad\lim_{|x|\to\infty} |x|^\alpha |f(x)|=\infty  
\quad\mbox{and}\quad 
\displaystyle\lim_{|x|\to\infty}f^\prime(x)=0.
\end{equation*}
For $\mu>0$ we introduce the bounded twice continuous differentiable function
\begin{equation*}
F_\mu(x):=\exp\Big(\mu\int_0^x f(y)dy\Big),\quad x\in\mathds{R}.
\end{equation*}
We shall apply Theorem 10.2.1 in \cite{SV}, p. 254, to the diffusion $X$, solution of (\ref{eqsign}), with the function $F_\mu$ for some $\mu>0$. It will implies that $\mathds{P}(\tau_e\leq T)>0$ for any $T>1$.  We need to verify that there exists $\lambda>0$ and $\mu>0$ such that for all $t\in[1,T]$ and $x\in \mathds R$,
\begin{equation}\label{minoration}
LF_\mu(t,x) \geq \lambda F_\mu(x)\quad\mbox{and}\quad \ln\Big(\frac{\sup_{x\in\mathds R} F_\mu(x)}{F_\mu(x_0)}\Big)<\lambda(T-1).
\end{equation}
Here $L$ is given in (\ref{ig}). In order to prove (\ref{minoration}), note that for all $t\in[1,T]$ and $x\in\mathds R$,
\begin{equation*}
LF_\mu(t,x)=\mu F_\mu(x)\Big(\rho t^{-\beta}|x|^\alpha |f(x)| + \frac{\mu}{2} f^2(x)+\frac{1}{2}f^\prime(x)\Big).
\end{equation*}
The assumptions on $f$ imply that there exists $r\geq 1$ such that, for all $\mu>0$,
\begin{equation*}
LF_\mu \geq \frac{k}{2}\,\mu F_\mu \quad\mbox{on}\quad [1,T]\times \left(\left[-1/2k,1/2k\right]\cup\left[-r,r\right]^c\right).
\end{equation*}
Besides, since $f^2$ is bounded away from zero, while $|f^\prime|$ is bounded on $\left[-1/2k,-r\right]\cup \left[1/2k,r\right]$, we deduce that there exists $\mu_0>0$ such that
\begin{equation*}
LF_{\mu_0} \geq \frac{k}{2}\,\mu_0 F_{\mu_0} \quad\mbox{on}\quad [1,T]\times \left(\left[-1/2k,-r\right]\cup\left[1/2k,r\right]\right).
\end{equation*}
Hence, for all $t\in[1,T]$ and $x\in\mathds R$, $LF_{\mu_0}(t,x)\geq \frac{k}{2}\,\mu_0 F_{\mu_0}(x)$ and we can see that
\begin{equation*}
 \ln\left(\frac{\sup_{x\in\mathds R} F_{\mu_0}(x)}{F_{\mu_0}(x_0)}\right)=\mu_0\int_{|x_0|}^\infty f(y)dy\leq \mu_0 \int_0^\infty g(y)dy<\frac{k}{2}\,\mu_0(T-1).
\end{equation*}
Therefore Theorem 10.2.1 in \cite{SV} 
%p. 254, 
applies with $\lambda:=\frac{k}{2}\mu_0$ and $F_{\mu_0}$ and $X$ explodes in finite time with positive probability. 
This ends the proof of the proposition, excepted for Lemma \ref{lemmepresquehorrible}.\hfill$\Box$

%%%%%%%%%%%%%%%%%%%%%%%%%%%%%%%%%%%%%%%%%%%%%%%%%%%%%%%%%%%%%%%%%%%%%%%%%%%%%%%%%%%%%%%%%%%%

\section{Asymptotic behaviour of solutions}
\setcounter{equation}{0}

We present here the systematic study of the recurrence, transience or convergence of the time-inhomogeneous one-dimensional diffusion $X$ (a regular strong Markov process solution of (\ref{eqsign}))
for parameters $(\rho,\alpha,\beta)\in\mathcal P:=\mathcal P_-\cup\mathcal P_+\,$,
where
\begin{equation*}
\mathcal P_{-} := (-\infty,0)\times (-1,\infty)\times \mathds R\;\;\mbox{(attractive case)}\;\;\mbox{and}\;\;\mathcal P_{+} :=  (0,\infty)\times \mathds R \times \mathds R\;\;\mbox{(repulsive case)}.
\end{equation*}
Set $E_\alpha:=\mathds R$, when $\alpha\in (-1,\infty)$, $E_\alpha:=(0,\infty)$, when $\alpha\in(-\infty,-1]$, and 
introduce the probability distributions $\Lambda_{\rho,\alpha}$ and $\Pi_{\rho,\alpha}$ on $E_\alpha$  defined by
\begin{equation}\label{Lambda}
\Lambda_{\rho,\alpha}(dx):=c^{-1}e^{-V_{\rho,\alpha}(x)}e^{-x^2/2}\,dx
\quad\mbox{and}\quad
\Pi_{\rho,\alpha}(dx):=k^{-1}e^{-V_{\rho,\alpha}(x)}\,dx.
\end{equation}
Here we denote $c,k$ the normalization constants and
\begin{equation}\label{simplepotential}
V_{\rho,\alpha}(x):=V_{\rho,\alpha,\beta}(1,x)=V_{\rho,\alpha,0}(t,x)
=\left\{\begin{array}{ccc}
-\frac{2\rho}{\alpha+1}|x|^{\alpha+1},&\mbox{ if }&\alpha\neq 1,\\
-2\rho\log{|x|},&\mbox{ if }&\alpha=1,\end{array}\right.
\end{equation}
where $V_{\rho,\alpha,\beta}(t,x)$ is the time-dependent potential given in (\ref{potential}). Besides, let us introduce the following three rate functions,
\begin{equation}\label{LLI1}
L(t):=(2t\ln{\ln{t}})^{\frac{1}{2}},\;\;
L_{\rho,\alpha}(t):=t^{\frac{1}{2}} \left(c_{\rho,\alpha}\ln{\ln{t}}\right)^{\frac{1}{\alpha+1}}\;\;\mbox{and}\;\; L_{\rho,\alpha,\beta}(t):=(c_{\rho,\alpha,\beta}\,t^\beta\ln{t})^{\frac{1}{\alpha+1}},
\end{equation}
where
\begin{equation}\label{constants}
\quad c_{\rho,\alpha}:=\frac{|\alpha+1|}{2|\rho|}\quad\mbox{and}\quad c_{\rho,\alpha,\beta}:=\frac{|\alpha+1-2\beta|}{2|\rho|}.
\end{equation}

\noindent
We shall say that the process $X$ is {\sl recurrent} in $E\subset \mathds R$ 
if, for all $x\in E$, the set $\{t\geq 1 : X_t = x \}$ is unbounded a.s. and we shall say that it is {\sl transient} 
if $\lim_{t\to \tau_e}|X_t|=\infty$ a.s.

%%%%%%%%%%%%%%%%%%%%%%%%%%%%%%%%%%%%%%%%%%%%%%%%%%%%%%%%%%%%%%%%%
 
\subsection{Behaviour on the critical line: $2\beta=\alpha+1$}

The scaling transformation (\ref{scalingt}) associated with the exponential 
change of time provides a time-homogeneous equation (\ref{criticalline}). With the help of Motoo's theorem 
(see Theorem \ref{Motoo} below) and of the ergodic theorem (see, for instance, Theorem 23.15 in \cite{Ka}, p. 465) we obtain the asymptotic behaviour of solutions 
to (\ref{eqsign}).

\begin{thm}[Attractive case]\label{criticallineattract}
If $(\rho,\alpha,\beta)\in \mathcal P_-$ and $2\beta=\alpha+1$, $X$ is recurrent in $\mathds R$ and 
\begin{equation}\label{criticlineattractlaw}
\lim_{t\to\infty}\frac{X_t}{\sqrt t}\overset{\mathcal L}{=}\Lambda_{\rho,\alpha}.
\end{equation}
Moreover,
\begin{enumerate}
 \item[i)] if $\alpha\in(-1,1)$, it satisfies
\begin{equation}\label{criticlineattractLLI11}
\limsup_{t\to\infty}\frac{X_t}{L(t)}=1\;\;\mbox{a.s.};
\end{equation}
 \item[ii)] if $\alpha\in(1,\infty)$, it satisfies
\begin{equation}\label{criticlineattractLLI12}
\limsup_{t\to\infty}\frac{X_t}{L_{\rho,\alpha}(t)}=1\;\;\mbox{a.s.};
\end{equation}
 \item[iii)] if $\alpha=1$, it satisfies
\begin{equation}\label{criticlineattractLLI2}
\limsup_{t\to\infty}\frac{X_t}{L(t)}=\frac{1}{\sqrt{1-2\rho}}\quad\mbox{a.s.}
\end{equation}
\end{enumerate}
\end{thm}

\noindent
In the repulsive case similar ideas will apply. However, we need to distinguish two particular 
cases, when $\alpha=-1$ (Bessel case) or when $\alpha=1$ (the continuous time analogue of the Friedman's 
urn model in \cite{Fre} and \cite{MV}). We note that in these cases the recurrent or transient features 
depend on the position of $\rho$ with respect to $1/2$.

\begin{thm}[Repulsive case]\label{criticallinerepuls}
Assume that $(\rho,\alpha,\beta)\in \mathcal P_+$  and $2\beta=\alpha+1$. 
\begin{enumerate}
 \item[i)] If $\alpha\in(-1,1)$, $X$ is recurrent in $\mathds{R}$ and it satisfies (\ref{criticlineattractlaw}) and
(\ref{criticlineattractLLI11}).
\item[ii)] If $\alpha\in(-\infty,-1)$, $X$ is transient, it satisfies (\ref{criticlineattractlaw}), (\ref{criticlineattractLLI11}) and 
\begin{equation}\label{criticlinerepulsLLItrans}
%~ \lim_{t\to\infty}\frac{X_t}{\sqrt t}\overset{\mathcal L}{=}\Lambda_{\rho,\alpha},\quad
%~ \limsup_{t\to\infty}\frac{X_t}{L(t)}=1\quad\mbox{a.s.}\quad\mbox{and}\quad
\liminf_{t\to\infty}\frac{X_t}{L_{\rho,\alpha}(t)}=1\quad\mbox{a.s.}. 
\end{equation}
\item[iii)] If $\alpha\in(1,\infty)$, the explosion time $\tau_e$ of $X$ is finite a.s. and 
\begin{equation}\label{vitesseexplosion}
|X_t|\underset{t\to\tau_e}{\sim}\frac{ \tau_e^{\frac{\alpha+1}{2(\alpha-1)}}}{(\rho(\alpha-1)(\tau_e-t))^{\frac{1}{\alpha-1}}}\quad\mbox{a.s.}
\end{equation}
\item[iv)] If $\alpha=-1$, $X$ is the classical Bessel process of dimension $2\rho+1$. It satisfies (\ref{criticlineattractlaw}) and (\ref{criticlineattractLLI11})  and, it is recurrent in $[0,\infty)$, when $\rho\in(0,1/2)$, recurrent in $(0,\infty)$, when $\rho=1/2$ and transient, when $\rho\in(1/2,\infty)$. 
Moreover, 
\begin{equation}\label{besseltrans}
\liminf_{t\to\infty}\frac{\ln\Big(\frac{X_t}{\sqrt t}\Big)}{\ln\ln t}= -\frac{1}{2\rho -1}\quad\mbox{a.s.}\quad\mbox{when}\quad\rho\in(1/2,\infty).
\end{equation}
\item[v)] If $\alpha=1$,  $X$ is a Gaussian process, recurrent in $\mathds R$, when $\rho\in(0,1/2]$, and transient, when $\rho\in(1/2,\infty)$. Moreover,
\begin{enumerate} 
\item[a)] if $\rho\in(0,1/2)$, it satisfies
\begin{equation*}\label{linearcase0}
\lim_{t\to\infty}\frac{X_t}{\sqrt{t}}\overset{\mathcal L}{=}\mathcal N\left(0,\frac{1}{2\rho-1}\right)\quad\mbox{and}
\quad\limsup_{t\to\infty}\frac{X_t}{L(t)}=\sqrt{\frac{2}{1-2\rho}}\quad\mbox{a.s.};
\end{equation*}
\item[b)] if $\rho=1/2$, it satisfies
\begin{equation*}\label{linearcase1}
\lim_{t\to\infty}\frac{X_t}{\sqrt{t \ln t}}\overset{\mathcal L}{=}\mathcal N(0,1)\quad\mbox{and}\quad \limsup_{t\to\infty}\frac{X_t}{\sqrt{2\,t\ln t \ln\ln\ln t}}=1\quad\mbox{a.s.};
\end{equation*}
\item[c)] if $\rho\in(1/2,\infty)$, it satisfies
\begin{equation*}\label{linearcase2}
\lim_{t\to\infty} \frac{X_t}{t^\rho}=G_{\rho,x_0}\quad\mbox{a.s.,}\quad\mbox{with}\;\; G_{\rho,x_0}\sim\mathcal N\left(x_0,\frac{1}{2\rho-1}\right).
\end{equation*}
\end{enumerate}
\end{enumerate}
\end{thm}

\begin{remarq}
The results contained in the latter theorem are in keeping with some results obtained for discrete time models in \cite{Fre} and \cite{MV}. More precisely, for the case $(\alpha,\beta)=(1,1)$ (point v) of Theorem \ref{criticallinerepuls})
one finds similar results as Theorems 3.1, 4.1 and 5.1 in \cite{Fre} and Corollary 1 in \cite{MV} concerning Friedman's urn model. For the case $(\alpha,\beta)=(-1,0)$ (point iv) of Theorem \ref{criticallinerepuls}) one gets similar result as in Theorem 5 in \cite{MV}. We also point out that the part i) of Theorem \ref{criticallinerepuls}  gives the asymptotic behaviour on the domain where the question is stated as an open problem in \cite{MV}, p. 958. 
\end{remarq}

\noindent
{\bf Proof of Theorem \ref{criticallineattract}.} Let $Z=\Phi_{\mathtt e}(X)\equiv \Phi_{\upgamma}(X)$ be the solution of the time-homogeneous equation (\ref{criticalline}). The scale function and the speed measure of $Z$ (see Chap. VI in \cite{IW}, pp. 446-449) are respectively given by
\begin{equation*}
s(x):=\int_0^{x}e^{V_{\rho,\alpha}(y)}e^{\frac{y^2}{2}}\,dy
\quad\mbox{and}\quad
m(dx):=e^{-V_{\rho,\alpha}(x)}e^{-\frac{x^2}{2}}\,dx.
\end{equation*}
Remark that $m$ is a finite measure on $\mathds{R}$ and $m(dx)/m(\mathds R)=\Lambda_{\rho,\alpha}(dx)$.  By using 
the ergodic theorem (see, for instance, Theorem 23.15 in \cite{Ka}, p. 465), we obtain
\begin{equation*}
\lim_{t\to\infty} \frac{X_t}{\sqrt t}\overset{}{=}
\lim_{t\to\infty} Z_{\ln{t}}\overset{\mathcal L}{=} \Lambda_{\rho,\alpha}.
\end{equation*}
To complete the proof we shall apply Motoo's theorem (see \cite{Mo}). We recall this result since it will be used several times.

\begin{thm}[Motoo]\label{Motoo}
Let $X$ be a regular continuous strong Markov process in $(a,\infty)$, $a\in[-\infty,\infty)$, 
which is homogeneous in time,  with scale function $s$ and finite speed measure $m$. For every real positive increasing function $h$,
\begin{equation*}
\mathds{P}\Big(\limsup_{t\to\infty} \frac{X_t}{h(t)}\geq 1\Big)=0\quad \mbox{or}\quad 1
\quad\mbox{ according to whether }\quad\int^{\infty}\frac{dt}{s(h(t))}<\infty\quad\mbox{or}\quad
=\infty.
\end{equation*}
\end{thm}

\noindent
Recall that $V_{\rho,\alpha}$ is given by (\ref{simplepotential}).
By using L'H{\^o}pital's rule, we can see that
\begin{equation*}
s(x)
\underset{x\to\infty}{\sim} \left\{
\begin{array}{ccc}
 x^{-1} e^{-V_{\rho,\alpha}(x)}e^{-x^2/2},&\mbox{ if }& \alpha\in(-1,1),\\
{(2|\rho| x^\alpha)}^{-1}e^{-V_{\rho,\alpha}(x)}e^{-x^2/2},&\mbox{ if } & \alpha\in(1,\infty).
\end{array}\right.
\end{equation*}

If $\alpha\in(-1,1)$, by a simple application of the Motoo's theorem we see that, for all $\ve>0$,
\begin{equation*}
\mathds{P}\Big(\limsup_{t\to\infty}\frac{Z_t}{\sqrt{2\ln t}}\geq 1+\ve\Big) =0\quad\mbox{and}\quad 
\mathds{P}\Big(\limsup_{t\to\infty}\frac{Z_t}{\sqrt{2\ln t}}\geq 1-\ve\Big)=1.
\end{equation*}
We deduce
\begin{equation*}
\limsup_{t\to\infty}\frac{X_t}{L(t)} = \limsup_{t\to\infty}\frac{Z_{\ln t}}{\sqrt{2\ln t}}=1\quad \mbox{a.s.}
\end{equation*}

If $\alpha\in(1,\infty)$, we deduce, again by Motoo's theorem,
\begin{equation*}
\limsup_{t\to\infty}\frac{X_t}{L_{\rho,\alpha}(t)} =  \limsup_{t\to\infty}\frac{Z_{\ln t}}{(c_{\rho,\alpha}\ln \ln t)^{\frac{1}{\alpha+1}}}=1\quad \mbox{a.s.}
\end{equation*}

Finally, assume that $\alpha=1$ (the linear case). Equality (\ref{criticlineattractLLI2}) can be proved by using similar methods as previously or by using standard results on linear stochastic differential equations.

Furthermore, by symmetry of equation (\ref{eqsign}), we can replace $X$ by $-X$ in relations 
(\ref{criticlineattractLLI11})-(\ref{criticlineattractLLI2}) to deduce that
$\limsup_{t\to\infty} X_t=\infty$ and $\liminf_{t\to\infty} X_t=-\infty$ a.s. and conclude that $X$ is recurrent in $\mathds R$.\hfill$\Box$\\

%%%%%%%%%%%%%%%%%%%%%%%%%%%%%%%%%%%%%%%%%%%%%%%%%%%%%%%%%%%%%%

\noindent
{\bf Proof of Theorem \ref{criticallinerepuls}.} To begin with, we point out that the proof of the point {\sl i)},
when $\alpha\in(-1,1)$, is the same as in the proof of Theorem \ref{criticallineattract}. 

When $\alpha\in(-\infty,-1)$, this last statement is also true when proving (\ref{criticlineattractlaw}) and 
(\ref{criticlineattractLLI11}). We need to prove (\ref{criticlinerepulsLLItrans}) and the transient feature.
To this end, consider again $Z=\Phi_{\mathtt e}(X)\equiv \Phi_{\upgamma}(X)$. By Ito's formula, we can see that $\tilde Z$ is the weak solution of
\begin{equation*}\label{liminf}
d{\tilde Z}_t= {\tilde Z}_t^2\, dW_t +\left({\tilde Z}_t^3 -\rho {\tilde Z}_t^{2-\alpha} + \frac{{\tilde Z}_t}{2}\right)\,dt,\quad {\tilde Z}_0=\frac{1}{x_0},\quad\mbox{with}\quad {\tilde Z}:=\frac{1}{Z}.
\end{equation*}
Again, by applying Motoo's theorem to $\tilde Z$, we deduce that 
\begin{equation*}
\liminf_{t\to\infty}\frac{X_t}{L_{\rho,\alpha}(t)}=\bigg(\limsup_{t\to\infty}\frac{{\tilde Z}_t}{(c_{\rho,\alpha}\ln t)^{\frac{1}{|\alpha +1|}}}\bigg)^{-1}
=1\quad\mbox{a.s.}
\end{equation*}
Note that this relation insure the transient feature, since $\lim_{t\to\infty}L_{\rho,\alpha}(t)=\infty$, when $\alpha<-1$.

Assume that $\alpha\in(1,\infty)$. We have already showed that the explosion time $\eta_e$ of $Z$ 
is finite a.s. (see Proposition \ref{timeexplosion}). Moreover, we can see that the process $z_t:=Z_t-W_t$ satisfies the random ordinary differential equation 
\begin{equation*}
\frac{dz_t}{dt} = \rho\,{\rm sgn}(z_t+W_t)|z_t+W_t|^{\alpha} -\frac{z_t+W_t}{2}.
\end{equation*}
We deduce that
\begin{equation*}
\frac{|z_t|^{1-\alpha}}{\alpha-1}=\int_t^{\eta_e}\frac{dz_s}{{\rm sgn}(z_s)|z_s|^\alpha}\underset{t\to\eta_e}{\sim} \rho (\eta_e-t)\quad\mbox{and}\quad |Z_t| \underset{t\to\eta_e}{\sim}\frac{1}{(\rho(\alpha-1)(\eta_e-t))^{\frac{1}{\alpha-1}}}\quad\mbox{a.s.}
\end{equation*}
Remark also that the explosion time $\tau_e$ of $X$ satisfies $\tau_e=e^{\eta_e}$ a.s. Therefore
\begin{equation*}
|X_t|=\sqrt{t}\,|Z_{\ln t}|\underset{t\to\tau_e}{\sim}\frac{\sqrt{\tau_e}}{(\rho(\alpha-1)(\ln\tau_e-\ln t))^{\frac{1}{\alpha-1}}}\underset{t\to\tau_e}{\sim}\frac{ \tau_e^{\frac{\alpha+1}{2(\alpha-1)}}}{(\rho(\alpha-1)(\tau_e-t))^{\frac{1}{\alpha-1}}}\quad\mbox{a.s.}
\end{equation*}

Assume that $\alpha=-1$ and let $R^{(\mathtt e)}$ be the pathwise unique strong solution of equation (\ref{eqtimechangebesselcarre}). By applying Lemma 2.2 in \cite{AW}, p. 916 %(using Motoo's Theorem) 
and the ergodic theorem to $R^{(\mathtt e)}$, we obtain (\ref{criticlineattractlaw}) and (\ref{criticlineattractLLI11}) by change of time. Equality (\ref{besseltrans}) is a consequence of Lemma 4.1 in \cite{AW}, p. 926. The recurrent or the transient features are proved in Chap. IX in \cite{RY}.

Finally, if $\alpha=1$ we are studying the classical case of a linear stochastic differential equation.
By standard arguments (see, for instance, \cite{RY} Proposition 2.3, Chap. IX, p. 378, and Theorem 1.7, Chap. V, p. 182) there exists a Brownian motion $W$ such that
\begin{equation*}
\frac{X_t}{t^{\rho}}=x_{0}+\int_{1}^{t}\frac{dB_s}{s^{\rho}}=x_{0}+W_{\phi(t)},
\quad\mbox{with}\quad 
\phi(t):=\left\{\begin{array}{ccc}
\frac{t^{1-2\rho}-1}{1-2\rho}&\mbox{if}& \rho\neq {1}/{2}\\
\ln{t}& \mbox{if} & \rho={1}/{2}.\end{array}\right. 
\end{equation*}
By using the well known properties of the Brownian motion, we deduce the convergence in distribution and the pathwise 
largest deviations of the Gaussian process $X$. Furthermore, the  recurrent or transient features are simple consequences.
\hfill$\Box$

%%%%%%%%%%%%%%%%%%%%%%%%%%%%%%%%%%%%%%%%%%%%%%%%%%%%%%%%%%%%%%%%%%

\subsection{Behaviour above the critical line: $2\beta>\alpha+1$}

The scaling transformation (\ref{scalingt}) associated with the exponential change of time does 
not provides a time-homogeneous equation. However, we shall prove that the asymptotic behaviour 
of equation (\ref{eqtimechangeexp}) is related to the asymptotic behaviour of the Ornstein-Uhlenbeck process 
(\ref{OU}), with the help of the Motoo theorem, the ergodic theorem, the comparison theorem (see, for instance, Theorem 1.1, Chap. VI in \cite{IW}, p. 437) and of the following result, whose proof is postponed to the Appendix.

\begin{lem}\label{convergence_distribution}
Let $Z$ and $H$ be regular strong Markov processes which are, respectively, weak solutions of the stochastic differential equations with continuous coefficients:
\begin{equation*}
dZ_s=\sigma(s,Z_s)\,dB_t+d(s,Z_s)\,ds\quad\mbox{and}\quad dH_s=\sigma_{\infty}(H_s)\,dB_s+d_{\infty}(H_s)\,ds.
\end{equation*}
Assume $(Z,H)$ is  asymptotically time-homogeneous and $\Pi$-ergodic, in the sense that
\begin{equation*}
\lim_{s\to\infty}\sigma(s,z)=\sigma_{\infty}(z)
\quad\mbox{and}\quad
\lim_{s\to\infty}d(s,z)=d_{\infty}(z),
\quad\mbox{uniformly on compact subsets of }\mathds{R},
\end{equation*} 
and $H$ converges in distribution to $\Pi$.
Furthermore, assume that $Z$ is bounded in probability, that is, for all $\varepsilon>0$ there exists $r>0$ such that 
$\sup_{s\geq 0} \mathds{P}(|Z_s|\geq r)<\varepsilon$.
Then $Z$ converges also in distribution to $\Pi$. 
\end{lem}

\begin{thm}[Attractive case]\label{abovecriticallineattract}
If $(\rho,\alpha,\beta)\in\mathcal P_-$ and $2\beta\in(\alpha+1,\infty)$, $X$ is recurrent in $\mathds R$ and
\begin{equation}\label{abovecriticattractlawLLI}
\lim_{t\to\infty}\frac{X_t}{\sqrt t}\overset{\mathcal L}{=}\mathcal N(0,1)\quad\mbox{and}\quad \limsup_{t\to\infty}\frac{X_t}{L(t)}=1\quad\mbox{a.s.}
\end{equation}
\end{thm}

\noindent
One more time, for the repulsive case we will follow similar ideas as for the attractive case, by modifying 
some computations when technical difficulties appear. 
Besides, when $\alpha\in(1,\infty)$ the process explodes with 
positive probability (see Proposition \ref{timeexplosion2}). Hence we need to adapt 
Lemma \ref{convergence_distribution} to show that, under the conditional probability of nonexplosion, the solution 
of equation (\ref{eqtimechangeexp}) behaves as the Ornstein-Uhlenbeck process (\ref{OU}).

\begin{thm}[Repulsive case]\label{abovecriticallinerepuls}
Assume that $(\rho,\alpha,\beta)\in\mathcal P_+$ and $2\beta\in(\alpha+1,\infty)$. 
\begin{enumerate}
\item[i)] If $\alpha\in(-1,1]$, $X$ is recurrent in $\mathds R$ and it satisfies (\ref{abovecriticattractlawLLI}).
\item[ii)] If $\alpha\in(-\infty,-1]$, $X$ is recurrent in $[0,\infty)$, when $\alpha=-1$, in $(0,\infty)$, when $\alpha\in(-\infty,-1)$ and $\beta\in[0,\infty)$ and it is transient, when $\alpha\in(-\infty,-1)$ and $\beta\in(-\infty,0)$. Moreover,
\begin{equation}\label{abovecriticallinerepulslaw}
\lim_{t\to\infty}\frac{X_t}{\sqrt t}\overset{\mathcal L}{=}|G|\quad\mbox{ and }\quad
\limsup_{t\to\infty}\frac{X_t}{L(t)}=1\quad\mbox{a.s.,}\quad{with}\quad G\sim\mathcal N(0,1),
\end{equation}
and
\begin{equation}\label{abovecriticallinerepulsLLI}
\liminf_{t\to\infty}\frac{X_t}{L_{\rho,\alpha,\beta}(t)}\geq 1\quad\mbox{a.s.}\quad\mbox{when}\quad \alpha\in(-\infty,-1).
\end{equation}
\item[iii)] If $\alpha\in(1,\infty)$,  conditionally to  $\{\tau_e=\infty\}$, $X$ is recurrent in $\mathds R$ and it satisfies (\ref{abovecriticattractlawLLI}), and conditionally to $\{\tau_e<\infty\}$, it satisfies  (\ref{vitesseexplosion}).
\end{enumerate}
\end{thm}

\begin{remarq}
The preceding statement concerning the recurrent asymptotic behaviour is a similar result as Theorem 4.2 ii) in \cite{MV}, p. 955. 
\end{remarq}

\noindent
{\bf Proof of Theorem \ref{abovecriticallineattract}.} The equalities in the statement 
will be consequences of Lemma \ref{convergence_distribution} and Motoo's theorem. 
Let us consider $X^{(\mathtt e)}=\Phi_{\mathtt e}(X)$ the unique weak solution of (\ref{eqtimechangeexp}) 
and $U$ the Ornstein-Uhlenbeck process solution of (\ref{OU}). Equalities in 
(\ref{abovecriticattractlawLLI}) are equivalent to
\begin{equation}\label{abovecriticattractlawLLI2}
\lim_{t\to\infty}X^{(\mathtt e)}_t\overset{\mathcal L}{=}G\quad\mbox{and}\quad \limsup_{t\to\infty}\frac{X^{(\mathtt e)}_t}{ \sqrt{2\ln t}}=1\quad\mbox{a.s.,}\quad\mbox{with}\quad G\sim\mathcal N(0,1).
\end{equation}
Equalities (\ref{abovecriticattractlawLLI2}) are satisfied by $U$ and roughly speaking $X^{(\mathtt e)}$ behaves 
as the Ornstein-Uhlenbeck process $U$. The proof of (\ref{abovecriticattractlawLLI2}) is split in three steps.

{\sl Step a).} Lemma \ref{convergence_distribution} does not %always 
apply to $(X^{(\mathtt e)},U)$ since the coefficients of  equation (\ref{eqtimechangeexp}) are discontinuous 
when $\alpha<0$. To remove the singularity, we consider $C:=(X^{(\mathtt e)})^3$ and $Q:=U^3$.
Ito's formula allows to see that $C$ and $Q$ are solutions of 
\begin{equation}\label{C}
dC_t=3\,C_t^{\frac{2}{3}}\,dW_t +3\Big(C_t^{\frac{1}{3}} -\frac{ C_t}{2} +\rho\,e^{\left(\frac{\alpha+1}{2}-\beta\right)t}\,{\rm sgn}(C_t)|C_t|^{\frac{\alpha+2}{3}}\Big)\,dt ,\quad C_0=x_0^3,
\end{equation}
and, respectively,
\begin{equation}\label{Q}
dQ_t=3\,Q_t^{\frac{2}{3}}\,dW_t +3\Big(Q_t^{\frac{1}{3}} -\frac{ Q_t}{2}\Big)\,dt,\quad Q_0=x_0^3. 
\end{equation} 
Since $2\beta>\alpha+1$, we deduce that $(C,Q)$ is asymptotically time-homogeneous and ${\mathcal L}(G^3)$-ergodic. 

{\sl Step b).} In order to apply Lemma \ref{convergence_distribution} to $(C,Q)$ we need to show that $C$ is bounded in probability. We prove this result by comparison with time-homogeneous ergodic diffusions. To this end, consider the  pathwise unique strong solution $C$ of equation (\ref{C}) and denote by $C^{\pm }$ the pathwise unique strong solutions of equations
\begin{equation*}
dC^{\pm }_t=3\,\left(C^{\pm}_t\right)^{\frac{2}{3}}\,dW_t +3\Big(\left(C^{\pm }_t\right)^{\frac{1}{3}} -\frac{ C^{\pm }_t}{2}\mp\rho\,|C_t^\pm|^{\frac{\alpha+2}{3}} \mathds{1}_{\{\mp C^\pm_t\geq 0\}}\Big)\,dt ,\quad C_0^\pm=x_0^3. 
\end{equation*}
By using a comparison theorem (see Theorem 1.1, Chap. VI in \cite{IW}, p. 437) we get, 
for all $t\geq 0$, $C^-_t\leq C_t\leq C^+_t$, a.s. Moreover, by computation of the speed measure 
as in the proof of Theorem \ref{criticallineattract}, we can see that $C^{\pm}$ are ergodic diffusions 
and therefore they are bounded in probability. By comparison, it is the same for $C$, and this fact implies 
the first equality in (\ref{abovecriticattractlawLLI2}).

{\sl Step c).} We get the pathwise largest deviations of $C$ by comparison with the 
time-homogeneous ergodic diffusion. By applying Motoo's theorem to $C^+$ (as in Theorem \ref{criticallineattract}), we obtain
\begin{equation}\label{pluspetit}
\limsup_{t\to\infty}\frac{C_t}{(2\ln t)^{\frac{3}{2}}}\leq
\limsup_{t\to\infty}\frac{C_t^+}{(2\ln t)^{\frac{3}{2}}}=1
 \quad\mbox{a.s.}
\end{equation}
To deduce the second equality in (\ref{abovecriticattractlawLLI2}), we need to prove the opposite inequality in (\ref{pluspetit}). We can see that the equality (\ref{pluspetit})  holds for $-C_t$, by symmetry of  (\ref{C}), and it implies that
\begin{equation}\label{negligible1}
\lim_{t\to\infty}\rho\,e^{\left(\frac{\alpha+1}{2}-\beta\right)t}\,{\rm sgn}(C_t)|C_t|^{\frac{\alpha+2}{3}}=0\quad\mbox{a.s.}
\end{equation}
Let $u\geq 0$ be and let us introduce the pathwise unique strong solution of equation
\begin{equation*}
dC_t(u)=3\,C_t(u)^{\frac{2}{3}}\,dW_t +3\Big(C_t(u)^{\frac{1}{3}} -\frac{ C_t(u)}{2} -1\Big)\,dt ,\quad C_u(u)=C_u. 
\end{equation*}
We shall prove that for all $t\geq u$,
\begin{equation}\label{comparison1}
C_{t}(u)\leq C_t\quad\mbox{a.s. on}\quad \Omega_u:=\Big\{\sup_{t\geq u}\,\rho\,e^{\left(\frac{\alpha+1}{2}-\beta\right)t} \,|C_t|^{\frac{\alpha+2}{3}}\leq 1\Big\}.
\end{equation}
Indeed, we introduce the stopping time $\tau_u$ defined by
\begin{equation*}
\tau_u:=\inf\Big\{t\geq u : \,\rho\,e^{\left(\frac{\alpha+1}{2}-\beta\right)t} \,|C_t|^{\frac{\alpha+1}{3}} >1\Big\}.
\end{equation*}
Using again the comparison theorem in \cite{IW}, p.437, and a classical argument of localisation, we obtain $C_{\Cdot\wedge \tau_u}(u)\leq C_{\Cdot\wedge \tau_u}$ a.s. Since $\{\tau_u=\infty\}=\Omega_u$ we deduce (\ref{comparison1}). By applying Motoo's theorem to $C(u)$ and by using (\ref{negligible1}), we get
\begin{equation*}
1=\limsup_{t\to\infty} \frac{C_t(u)}{(2\ln t)^{\frac{3}{2}}}\leq \limsup_{t\to\infty}\frac{C_t}{(2\ln t)^{\frac{3}{2}}}\quad\mbox{a.s. on}\quad \Omega_u,\quad\mbox{and}\quad \mathds P(\cup_{u\geq 0}\Omega_u)=1.
\end{equation*} 
The opposite inequality in (\ref{pluspetit}) is obtained and the proof of (\ref{abovecriticattractlawLLI2}) is finished.

Finally, to conclude that $X$ is recurrent in $\mathds R$, it suffices to replace $X$ by $-X$ in the second 
equality in (\ref{abovecriticattractlawLLI}). This is possible by symmetry of equation (\ref{eqsign}).\hfill$\Box$\\

\noindent
{\bf Proof of Theorem \ref{abovecriticallinerepuls}.} 
Let us note that, for $\alpha\in(-1,1]$, the proof is exactly the same as the proof of Theorem \ref{abovecriticallineattract}, 
while the proof for $\alpha\in(-\infty,-1]$ follows some similar steps. 

Let us only bring out the differences for $\alpha\in(-\infty,-1]$. 
We consider again the couple $(X^{(\mathtt e)},U)$ and we  perform 
the {\sl Step a)} in the proof of Theorem \ref{abovecriticallineattract} with $((X^{(\mathtt e)})^3, U^3)$ replaced  by $(|X^{(\mathtt e)}|^{|\alpha|+1}, 
|U|^{|\alpha|+1})=:(C,Q)$ (to avoid the singularity) . It follows that $(C,Q)$ is asymptotically time-homogeneous and ${\mathcal L}(|G|^{|\alpha|+1})$-ergodic. Here $Q$ and $C$ are weak solutions of 
\begin{equation*}
dQ_t=(|\alpha|+1)\,Q_t^{\frac{|\alpha|}{|\alpha|+1}}\,dW_t +\frac{|\alpha|+1}{2}\Big(|\alpha|\,Q_t^{\frac{|\alpha|-1}{|\alpha|+1}} - Q_t\Big)\,dt,\quad Q_0=x_0^{|\alpha|+1},
\end{equation*} 
and
\begin{equation*}
dC_t=(|\alpha|+1) C_t^{\frac{|\alpha|}{|\alpha|+1}}dW_t +\frac{|\alpha|+1}{2}\Big[|\alpha|C_t^{\frac{|\alpha|-1}{|\alpha|+1}} - C_t + 2\rho\,e^{\left(\frac{\alpha+1}{2}-\beta\right)t} C_t^{\frac{|\alpha|+\alpha}{|\alpha|+1}}\Big]\,dt,\,\, C_0=x_0^{|\alpha|+1}.
\end{equation*}
As in {\sl Step b)} in the proof of Theorem \ref{abovecriticallineattract} we can show that $C_t$ is bounded in probability, by comparing $C_t$ with the ergodic nonnegative diffusion satisfying
\begin{equation*}
dC_t^+=(|\alpha|+1) (C_t^+)^{\frac{|\alpha|}{|\alpha|+1}}dW_t +\frac{|\alpha|+1}{2}\Big(|\alpha|(C_t^+)^{\frac{|\alpha|-1}{|\alpha|+1}} - C_t^+ + 2\rho\,(C_t^+)^{\frac{|\alpha|+\alpha}{|\alpha|+1}}\Big)\,dt,\, C_0^+=x_0^{|\alpha|+1}.
\end{equation*}
Lemma \ref{convergence_distribution} applies and we get the first equality in (\ref{abovecriticallinerepulslaw}).
Finally, as in {\sl Step c)} of the cited proof, by applying Motoo's theorem to $Q_t$ and $C^+_t$ and by comparison theorem we can obtain 
\begin{equation*}
1=\limsup_{t\to\infty}\frac{Q_t}{(2\ln t)^{\frac{|\alpha|+1}{2} }}\leq \limsup_{t\to\infty}\frac{C_t}{(2\ln t)^{\frac{|\alpha|+1}{2} }}\leq \limsup_{t\to\infty}\frac{C_t^+}{(2\ln t)^{\frac{|\alpha|+1}{2} }}=1\quad\mbox{a.s.}
\end{equation*}
We deduce the second equality in (\ref{abovecriticallinerepulslaw}).

If $\alpha< -1$, let $\tilde X^{(\upgamma)}:=1/X^{(\upgamma)}$
be the pathwise unique nonnegative strong solution of
\begin{equation*}
d \tilde X^{(\upgamma)}_t= (\tilde X^{(\upgamma)}_t)^2\,dW_t+\bigg( (\tilde X^{(\upgamma)}_t)^3 -\rho\, (\tilde X^{(\upgamma)}_t)^{2-\alpha} + \frac{\gamma\, \tilde X^{(\upgamma)}_t }{2(1-(1-\gamma)t)}\bigg) dt,\quad \tilde X^{(\upgamma)}_0:=\frac{1}{x_0}.
\end{equation*}
Recall that $X^{(\upgamma)}$ is the pathwise unique nonnegative strong solution of (\ref{eqtimechangegamma}).  Consider also the pathwise unique nonnegative strong solution $\tilde Y$ of 
\begin{equation*}
d \tilde Y_t= \tilde Y_t^2\,dW_t+\left( \tilde Y_t^3 -\rho\, \tilde Y_t^{2-\alpha} + \frac{|\gamma|\, \tilde Y_t }{2}\right) dt,\quad \tilde Y_0:=\frac{1}{x_0}.
\end{equation*}
By comparison between $\tilde X^{(\upgamma)}$ and $\tilde Y$, and by applying Motoo's theorem to $\tilde Y$, we deduce
\begin{equation*}
\liminf_{t\to\infty}\frac{X_t}{L_{\rho,\alpha,\beta}(t)}=\left(\limsup_{s\to\infty}\frac{X^{(\upgamma)}_s}{(c_{\rho,\alpha}\ln s )^{\frac{1}{|\alpha+1|}}}\right)^{-1}\geq \left(\limsup_{s\to\infty}\frac{\tilde Y_s}{(c_{\rho,\alpha}\ln s )^{\frac{1}{|\alpha+1|}}}\right)^{-1}=1\quad\mbox{a.s.}
\end{equation*}
In previous relation the first equality was obtained by using the change of time $s=\varphi_\upgamma^{-1}(t)$ defined in (\ref{gamma}).
Moreover, if $\alpha=-1$, the point $0$ is recurrent for $X$. By (\ref{abovecriticallinerepulslaw}), we get the recurrent feature in $[0,\infty)$. Besides, we obtain from (\ref{abovecriticallinerepulsLLI}) that $X$ is transient, when $\beta\in(-\infty,0)$. 
Furthermore, if $\beta=0$, $X$ is an homogeneous diffusion and by standard criteria, using the scale function, we can see that $X$ is recurrent in $(0,\infty)$. If $\beta\in[0,\infty)$, by comparison theorem with 
the process obtained for $\beta=0$, we get that $X$ is recurrent in $(0, \infty)$. 
The proof of ii) is complete.

If $\alpha\in(1,\infty)$, consider $X^{(\upgamma)}=\Phi_{\upgamma}(X)$ and $b$ the respective 
solutions of equations (\ref{eqtimechangegamma2}) and (\ref{deltapontmansuy}). 
Denote by $\eta_e$ the explosion time of $X^{(\upgamma)}$ and recall that $\{\eta_e\geq t_1\}=\{\tau_e=\infty\}$, $\mathds{P}(\eta_e=t_1)=0$ (Lemma \ref{lemmepresquehorrible}) and  $\lim_{t\to t_1} b_{t}=0$ a.s. By using (\ref{girsa}), 
\begin{equation*}
\mathds P\Big(\lim_{t\to t_1} X^{(\upgamma)}_t=0,\, \eta_e\geq t_1\Big)=\mathds E\Big({\mathds 1}_{\big\{\displaystyle\lim_{t\to t_1} b_t=0\big\}}\,\mathcal E(t_1)\Big)
=\mathds E(\mathcal E(t_1))=\mathds P(\eta_e\geq t_1).
\end{equation*}
By change of time, we get
\begin{equation}\label{pontlimitnonexplosion}
\lim_{t\to\infty}\frac{X_t}{t^{\frac{\beta}{\alpha+1}}}=\lim_{t\to t_1} X^{(\upgamma)}_t=0\quad\mbox{a.s.}\quad\mbox{on}\quad \{\tau_e=\infty\}.
\end{equation}
Therefore
\begin{equation*}
\lim_{t\to\infty}\frac{1}{\sqrt t}\int_1^t\rho\,\frac{|X_s|^\alpha}{s^\beta}ds=\lim_{t\to\infty}\frac{1}{\sqrt t}\int_1^t\rho\,\left|\frac{X_s}{s^\frac{\beta}{\alpha+1}}\right|^\alpha s^{-\frac{ \beta}{\alpha+1}}\,ds=0\quad\mbox{a.s.}\quad\mbox{on}\quad \{\tau_e=\infty\}.
\end{equation*}
We deduce that $X$ satisfies the iterated logarithm law in (\ref{abovecriticattractlawLLI}) under the conditional probability of nonexplosion. Hence it is recurrent in $\mathds R$. We shall prove the convergence in distribution 
(\ref{abovecriticattractlawLLI}) under the conditional probability of nonexplosion. For this end, it suffices to show that
\begin{equation}\label{conv_conditionnelle}
\lim_{s\to\infty}\mathds{P}(X^{(\mathtt e)}_s>x\mid \sigma_{e}=\infty)
=\frac{1}{\sqrt{2\pi}}\int_{x}^{\infty}\exp\left(-\frac{y^{2}}{2}\right)dy.
\end{equation}
Here $\sigma_e$ denotes the explosion time of $X^{(\mathtt e)}=\Phi_{\mathtt e}(X)$, the solution of (\ref{eqtimechangeexp}). Note that Lemma \ref{convergence_distribution} does not apply directly to $(X^{(\mathtt e)},U)$, since $\sigma_e$ could be finite with positive probability. By using (\ref{pontlimitnonexplosion}), we remark that
\begin{equation}\label{comparisonfin}
\lim_{s\to\infty}\rho\,e^{\left(\frac{\alpha+1}{2}-\beta\right)s}|\Um_s|^\alpha=\lim_{t\to\infty}\rho\,t^{\left(\frac{1}{2}-\frac{\beta}{\alpha+1}\right)}\left|\frac{X_t}{t^\frac{\beta}{\alpha+1}}\right|^{\alpha}=0
\quad \mbox{a.s. on }\quad\{\sigma_e=\infty\}.
\end{equation}
Let $\varepsilon>0, v\geq0$ be and denote $U^{(\pm\ve)}$ the pathwise unique strong solutions of equations 
\begin{equation*}
dU_s^{(\pm\ve )} = dW_s -\frac{U_s^{(\pm \ve)}}{2}  ds \,\pm \varepsilon\,ds,\quad U_v^{(\pm\ve)}=\Um_{v}\mathds{1}_{\{\sigma_{e}>v\}}.
\end{equation*}
It is classical that $U^{(\pm\ve)}$ is Feller and ergodic. Furthermore, the strong mixing property holds (see \cite{Ka}, Theorem 20.20, p. 408), hence we obtain
\begin{equation}\label{melange}
\lim_{s\to\infty}\mathds{P}(U_s^{(\pm\ve)}>x\mid \Omega_v^{\ve})=F_{\mp\ve}(x):=\frac{1}{
\sqrt{2\pi}}\int_{x}^{\infty}\exp\left(-\frac{(y\mp\,\ve)^{2}}{2}\right)dy,
\end{equation}
with
\begin{equation*}
\Omega_v^{\ve}:=\Big\{\sup_{s\geq v}\;\rho\,e^{\left(\frac{\alpha+1}{2}-\beta\right)s}|\Um_s|^\alpha \leq\varepsilon\Big\}.
\end{equation*}
Similarly as for (\ref{comparison1}), we can show, by using comparison theorem and a classical argument of localisation, that, for all $s\geq v$,
$U_s^{(-\ve)}\leq \Um_s\leq U_s^{(+\ve)}$ a.s. on $\Omega_{v}^{\ve}$.
We deduce, from (\ref{melange}),
\begin{equation}\label{convconditionelle}
F_{+\ve}(x)
\leq\liminf_{s\to\infty}\mathds{P}(\Um_s>x\mid \Omega_v^{\ve})
\leq\limsup_{s\to\infty}\mathds{P}(\Um_s>x\mid \Omega_v^{\ve})\leq 
F_{-\ve}(x).
\end{equation}
Thanks to (\ref{comparisonfin}) the set of nonexplosion is $\{\sigma_{e}=\infty\}=\cup_{v\geq 0}\Omega_v^{\ve}$. 
Letting $v\to\infty$, and then $\ve\to 0$ in (\ref{convconditionelle}), we  deduce (\ref{conv_conditionnelle}).

To finish the proof, we need to study the process $X$ conditionally to $\{\tau_e<\infty\}$ and prove that it 
 satisfies (\ref{vitesseexplosion}). The method is the same as in the 
proof of Theorem \ref{criticallinerepuls}: we show that 
\begin{equation*}
|X_t^{(\mathtt e)}|\underset{t \to\infty}{\sim}\frac{1}{(\rho(\alpha-1)(\eta_e-t) )^{\frac{1}{\alpha-1}}}\quad\mbox{a.s. on}\quad \{\eta_e<\infty\},
\end{equation*}
and we conclude by change of time.\hfill$\Box$

\subsection{Behaviour under the critical line: $2\beta<\alpha+1$}

In the attractive case, by using similar techniques as in the proofs of theorems \ref{abovecriticallineattract} and \ref{abovecriticallinerepuls}
we shall prove that the asymptotic behaviour of equation (\ref{eqtimechangegamma3}) is related to the asymptotic behaviour of the time-homogeneous equation (\ref{undercriticallinehomogeneous}). By change of time, we shall obtain the asymptotic behaviour for (\ref{eqsign}).\\

\begin{thm}[Attractive case]\label{undercriticallineattract}
If $(\rho,\alpha,\beta)\in\mathcal P_-$ and $2\beta\in(-\infty,\alpha+1)$, $X$ is recurrent in $\mathds R$, if $\beta\in[0,\infty)$, and $X$ converges a.s. towards $0$, if $\beta\in(-\infty,0)$. Moreover,
\begin{equation}\label{undercriticallinelawLLI}
\lim_{t\to\infty}\frac{X_t}{t^{\frac{\beta}{\alpha+1}}}\overset{\mathcal L}{=}\Pi_{\rho,\alpha}\quad\mbox{and}\quad \limsup_{t\to\infty} \frac{X_t}{L_{\rho,\alpha,\beta}(t)}=1\quad\mbox{a.s.}
\end{equation}
\end{thm}

\begin{thm}[Repulsive case]\label{transience}
If $(\rho,\alpha,\beta)\in\mathcal P_+$ and $2\beta\in(-\infty,\alpha+1)$, $X$ is transient. Moreover,  
\begin{enumerate} 
\item[i)]if $\alpha\in(-\infty,1)$, it satisfies
\begin{equation}\label{transient_scale}
\lim_{t\to\infty}\frac{|X_t|}{t^\frac{1-\beta}{1-\alpha}}=\Big(\frac{\rho(1-\alpha)}{1-\beta}\Big)^{\frac{1}{1-\alpha}}\quad\mbox{a.s.};
\end{equation}
\item[ii)] if $\alpha\in(1,\infty)$, it satisfies
\begin{equation}\label{vitesseexplosion2}
|X_t|\underset{t \to\infty}{\sim}\frac{\varphi_{\upgamma}^\frac{\gamma}{\alpha-1}\circ\varphi_\upgamma^{-1}(\tau_e)\cdot\tau_e^{\frac{\gamma}{2}}}{(\rho(\alpha-1)(\tau_e-t) )^{\frac{1}{\alpha-1}}}\quad\mbox{a.s.,}
\end{equation}
where $\varphi_\upgamma$ and $\gamma$ are given in (\ref{gamma});
\item[iii)] if $\alpha=1$, it satisfies
\begin{equation}\label{translinear}
\lim_{t\to\infty}\frac{X_t}{\exp\left(\frac{\rho\,t^{1-\beta}}{1-\beta}\right)}=G\quad \mbox{ a.s.,}
\end{equation}
where $G\sim\mathcal{N}(m,\sigma^2)$, with $m:=x_0\exp(\frac{\rho}{\beta-1})$ and 
$\sigma^2:=\int_{1}^\infty \exp(\frac{2\rho \,s^{1-\beta}}{\beta-1})ds$.
\end{enumerate}
\end{thm}

\begin{remarq}
Again, one finds a similar result as in Theorem 1 ii) from \cite{MV}, p. 951, concerning the transient feature
of the process. 
\end{remarq}

\noindent
{\bf Proof of Theorem \ref{undercriticallineattract}.} Let $\Hm=\Phi_{\upgamma}(X)$ and $H$ be the solutions 
respectively of (\ref{eqtimechangegamma3}) and (\ref{undercriticallinehomogeneous}). Equalities in (\ref{undercriticallinelawLLI}) are equivalent to
\begin{equation}\label{equivalent}
\lim_{t\to\infty}\Hm_t\overset{\mathcal L}{=}S\quad\mbox{and}\quad \limsup_{t\to\infty}\frac{\Hm_t}{(c_{\rho,\alpha}\ln t)^{\frac{1}{\alpha+1}}}=1\quad\mbox{a.s.,}\quad\mbox{with}\quad S\sim \Pi_{\rho,\alpha}.
\end{equation} 
Note that $H$ satisfies these equalities. To prove (\ref{equivalent}) we can follow similar Steps {\sl a)-c)} as in the proof of Theorem \ref{abovecriticallineattract}, 
by considering $C:=(\Hm)^3$ and $Q:=H^3$, which are the pathwise unique strong solutions of
\begin{equation*}
dQ_t= 3\,Q_t^{\frac{2}{3}}\,dW_t +3\Big(\rho\,{\rm sgn}(Q_t)|Q_t|^{\frac{\alpha+2}{3}} + Q^{\frac{1}{3}}_t\Big) dt,\quad Q_0=x_0^3,
\end{equation*}
and
\begin{equation*}
dC_t= 3\,C_t^{\frac{2}{3}}\,dW_t +3\Big(\rho\,{\rm sgn}(C_t)|C_t|^{\frac{\alpha+2}{3}} + C^{\frac{1}{3}}_t -\frac{\gamma \, C_t}{2(1-(1-\gamma)t)} \Big) dt,\quad C_0=x_0^3.
\end{equation*}
As in {\sl Step a)} we get that $(C,Q)$ is asymptotically homogeneous and ${\mathcal L}(S^3)$-ergodic. The arguments of the corresponding {\sl Step b)} are going on as follows.
Since $\rho$ is negative, it is not difficult to prove that $X^2_t\leq{\tilde W}^2_t$, $t\geq 0$, where ${\tilde W}_t$ is a Brownian motion. We obtain by using the change of time $s=\varphi_\upgamma^{-1}(t)$ and the iterated logarithm law that  
\begin{equation}\label{thetagamma}
\limsup_{t\to\infty}\frac{C_t}{(2t\ln\ln t)^{\frac{3}{2}}}\overset{}{=} \limsup_{s\to\infty}\frac{X_s^3}{\left(\frac{2}{1-\gamma}\,s\ln\ln s\right)^{\frac{3}{2}}}\leq (1-\gamma)^3\quad\mbox{a.s.}
\end{equation}
Let $\theta\in(1/3,(\alpha+2)/3)$ be. Thanks to (\ref{thetagamma}), we deduce 
\begin{equation}\label{comparisonfin2}
\lim_{t\to\infty}  \frac{1}{|C_t|^\theta}\cdot\frac{\gamma \,C_t}{1+(1-\gamma)t}=0\quad\mbox{a.s.}
\end{equation} 
Let $v\geq 0$ be and introduce $C^{\pm}$ the pathwise unique strong solution of 
\begin{equation*}
dC_t^{\pm}= 3\,(C_t^\pm)^{\frac{2}{3}} dW_t +3\Big(\rho\,{\rm sgn}(C_t^\pm)|C_t^\pm|^{\frac{\alpha+2}{3}} + (C^\pm_t)^{\frac{1}{3}} \pm |C_t^\pm|^{\theta} \Big) dt,\quad C_v^\pm=C_v.
\end{equation*}
As for the proof of comparison (\ref{comparison1}) we can prove that for all $s\geq v$,
\begin{equation}\label{doubleinegalite}
C^-_t\leq C_t\leq C^+_t\quad\mbox{a.s. on}\quad 
\Omega_v:=\Big\{\sup_{t\geq v}\;\frac{1}{|C_t|^\theta}\cdot\frac{\gamma\, C_t}{1+(1-\gamma)t}\leq 1\Big\}.
\end{equation} 
By (\ref{comparisonfin2}), for any $\varepsilon>0$, we can choose $v\geq 0$ such that $\mathds{P}(\Omega_v)\geq1-\varepsilon$. Moreover, there exists $r\geq 0$ such that for all $t\geq v$, 
$\mathds{P}(|C^{\pm}_t|\geq r)\leq \varepsilon$  since $C^\pm$ is an ergodic diffusion (by computation of the speed measure). 
Combining the latter inequality with (\ref{doubleinegalite}) which holds on $\Omega_v$, we obtain that 
$\mathds{P}(|C_t|\geq r)\leq 2\ve$, for all $t\geq v$ and therefore we conclude that $C$ is bounded in probability. Finally, {\sl Step c)} is  a consequence of Motoo's theorem applied to $C^\pm$ and to the preceding comparison, 
\begin{equation*}
1=\limsup_{t\to\infty}\frac{C_t^-}{(c_{\rho,\alpha}\ln t)^{\frac{3}{\alpha+1}}}\leq \limsup_{t\to\infty}\frac{C_t}{(c_{\rho,\alpha}\ln t)^{\frac{3}{\alpha+1}}}\leq \limsup_{t\to\infty}\frac{C_t^+}{(c_{\rho,\alpha}\ln t)^{\frac{3}{\alpha+1}}}
= 1\quad\mbox{ a.s.}
\end{equation*}
This ends the proof of (\ref{equivalent}). To get the recurrence feature or the convergence toward 0 we use 
the second equality in (\ref{undercriticallinelawLLI}) with $X$ and $-X$.\hfill$\Box$\\

\noindent
{\bf Proof of Theorem \ref{transience}.} Assume that $\alpha\in(-\infty,1)$. To simplify the computations, let us denote the limit and the exponent of $t$ in 
(\ref{transient_scale}), respectively by
\begin{equation*}
\ell:=\Big(\frac{\rho(1-\alpha)}{1-\beta}\Big)^\frac{1}{1-\alpha}\quad\mbox{and}\quad
\nu:=\frac{1-\beta}{1- \alpha}\,.
\end{equation*}
If we set $S_t:=X_{t}^2/{t^{2\nu}}$, it suffices to verify  that 
$\lim_{t\to\infty}S_{t}=\ell^{2}$ a.s., that is, for all $\ve>0$,
\begin{equation}\label{but}
\limsup_{t\to\infty}\,S_t \leq \ell^{2} + 3\varepsilon 
\quad\mbox{ and }\quad
\liminf_{t\to\infty}\,S_t \geq \ell^{2} -3\varepsilon 
\quad \mbox{ a.s.}
\end{equation}
We shall prove only the first inequality in (\ref{but}), the second one being obtained in a similar way. 
We split the proof of this inequality in four steps.

{\sl Step a).} We begin by proving that, for all $\ve>0$, 
\begin{equation}\label{stepa}
\{t\geq 1:S_{t}\leq \ell^{2}+\ve\}\quad \mbox{is unbounded a.s.}
\end{equation}
 For this end, set $\eta_u:=\inf\left\{v\geq u : S_v\leq \ell^2+\varepsilon\right\}$, $u\geq 1$. Then, it suffices to prove that  for all $u\geq 1$ large enough, $\eta_u<\infty$ a.s. By using Ito's formula, we can see that 
\begin{equation}\label{ito_positive}
S_{t\wedge\eta_{u}} =  S_{u}+\int_{u}^{t\wedge\eta_{u}}LG(s,X_{s})ds
+\int_{u}^{t\wedge\eta_{u}}\partial_x G(s,X_s)\,dB_{s}
 :=  S_{u}+M_t+A_t,
\end{equation}
where $G(t,x):=x^{2}/{t^{2\nu}}$ and where $L$ is given by (\ref{ig}). Moreover, we can see that there exist $s_{0}\geq 1$ and $c>0$ such that, 
for all $s\geq s_{0}$ and $x\in\mathds{R}$, for which $G(s,x)\geq \ell^2+\ve$ and 
\begin{equation}\label{estim_LV}
LG(s,x)= \frac{2\rho }{s}
\left(G(s,x)^{\frac{\alpha-1}{2}}-\ell^{\alpha-1}\right)G(s,x)+\frac{1}{s^{2\nu}}\leq -\frac{c}{s}\leq 0. 
\end{equation}
This implies that the local martingale part $S_{u}+M$ of the non-negative semimartingale in (\ref{ito_positive}), together with $S_{\bullet\wedge\eta_{u}}$  itself, 
are nonnegative supermartingales for all $u\geq s_0$. Therefore, the bounded variation part $A$ will be a convergent process as the difference of two convergent supermartingales. Thanks to (\ref{estim_LV}), this is possible if and only if $\eta_u<\infty$ a.s.

{\sl Step b).} 
We introduce an increasing sequence of stopping times as follows: 
\begin{equation*}
\tau_1:=\inf\{t\geq s_{0} : 
S_t=\ell^2+{2\varepsilon}\},\quad 
\sigma_1:=\inf\left\{t\geq \tau_1 : 
S_t\in \left\{\ell^2+\varepsilon,\ell^2+3\varepsilon\right\}\right\}
\end{equation*}
and for every integer $n\geq 2$,
\begin{equation*}
\tau_{n}:=\inf\{t\geq \sigma_{n-1} : 
S_t=\ell^2+{2\varepsilon}\},\quad 
\sigma_{n}:=\inf\left\{t\geq \tau_n : 
S_t\in\left\{\ell^2+\varepsilon,\ell^2+3\varepsilon\right\}\right\}.
\end{equation*}

\begin{figure}[!ht]\label{stoppingtime}
 \centering
 \includegraphics[width=13cm,height=3.3cm]
{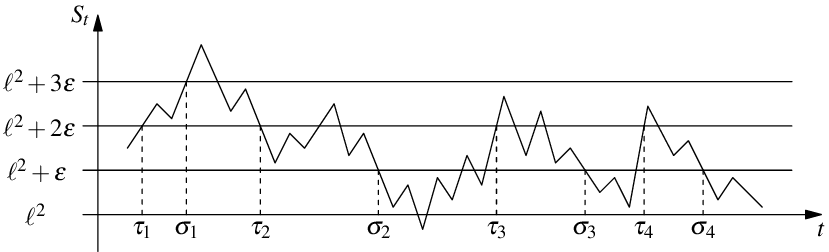}
%{sousmartingale-1.ps}
 \caption{The increasing sequence of stopping times}
 \end{figure}
 
\noindent
Set $F:=\cap_{n\geq 1}\{\tau_n<\infty\}$. Thanks to (\ref{stepa}), we obtain that
$\limsup_{t\to\infty}\,S_t \leq \ell^{2} + 2\varepsilon$ a.s. on $F^c$.
To prove the first inequality in (\ref{but}), we need to show that 
\begin{equation*}
\limsup_{t\to\infty}\,S_t \leq \ell^{2} + 3\varepsilon\quad\mbox{a.s. on }F,\quad\mbox{or, equivalently,}\quad \mathds{1}_{F}\sum_{n\geq 1}
\mathds{1}_{\{S_{\sigma_{n}}=\ell^{2}+3\ve\}}<\infty\quad\mbox{a.s.}
\end{equation*}
By using a conditional version of the Borel-Cantelli lemma (see, for instance, Corollary 7.20
in \cite{Ka}, p. 131), it is equivalent to prove that
\begin{equation}\label{but2}
\sum_{n=1}^\infty \mathds{P}\left(S_{{\sigma_n}}=\ell^2+3\ve,\tau_n<\infty\mid
\mathcal{F}_{\tau_n}\right)<\infty\quad\mbox{a.s.}
\end{equation} 

{\sl Step c).} We show that there exist positive constants $\lambda_1$ and $\lambda_2$ such that for all $n\geq 1$,
\begin{equation}\label{maj}
\mathds{P}\left(S_{{\sigma_n}}=\ell^2+3\ve,\tau_n<\infty\mid
\mathcal{F}_{\tau_n}\right)\leq \lambda_1 \tau_n^{\left(\frac{1}{2}-\nu\right)}\exp{\left(-\lambda_2 \tau_n^{2\nu-1}\right)}\quad \mbox{a.s. on }\quad \{\tau_n<\infty\}.
\end{equation}
To this end, let us denote by $\mathds{P}_{s,x}$ the distribution of the weak solution of (\ref{eqsign}) such that $X_s=x$.  The strong Markov property applies and this yields
\begin{equation*}
\mathds{P}\left(S_{\sigma_n}=\ell^2+3\ve,\tau_n<\infty\mid\mathcal{F}_{\tau_n}\right) 
 = \mathds{P}_{\tau_{n},X_{\tau_{n}}}\left( S_{\sigma_n}=\ell^2+3\ve\right)\quad\mbox{a.s. on}\quad \{\tau_n<\infty\}.
\end{equation*}
As in (\ref{ito_positive}), we can write under the conditional probability $\mathds{P}_{\tau_n,X_{\tau_n}}$, the canonical decomposition $S_{t\wedge\sigma_n} =  S_{\tau_n}+M_t^n +A^n_{t}$, sum of a local martingale and a bounded variation process.
Besides, we can show that $v_n:=\langle M^n\rangle_\infty$ satisfies
\begin{equation*}
v_n=\int_{\tau_n}^{\sigma_n}\frac{4 S_u}{u^{2\nu}}\,du \leq \int_{\tau_n}^\infty\frac{4(\ell^2+3\ve)}{u^{2\nu}}\,du=\frac{4(\ell^2+3\varepsilon)}{(2\nu -1)\tau_n^{2\nu-1}}\leq \frac{4(\ell^2+3\varepsilon)}{(2\nu -1)s_0^{2\nu-1}}=:v_0.
\end{equation*}
Then, by the Dambis-Dubins-Schwarz theorem, there exists a standard Brownian motion $W^n$ (under the conditional probability $\mathds{P}_{\tau_n,X_{\tau_n}}$) such that $M^n=W^n_{\langle M^n\rangle}$ and
since $A^n$ is strictly negative, we can see that 
\begin{equation*}
\Big\{\sup_{0\leq t\leq  v_n} W^n_t<\ve \Big\}
\subset
\Big\{\sup_{t\geq\tau_n} M_t^n<\ve\Big\}
\subset \left\{S_{\sigma_n}=\ell^2+\ve \right\}.
\end{equation*}
It is classical that
$\sup\{ W_t^n :0\leq t\leq v_n\} \overset{\mathcal L}=|G_n|$, with $G_n\overset{\mathcal L}{=}\mathcal{N}(0,v_n)$,
under the conditional probability $\mathds{P}_{\tau_n,X_{\tau_n}}$ and therefore we obtain 
\begin{equation*}
\mathds{P}_{\tau_n,X_{\tau_n}}(S_{\sigma_n}=\ell^2+3\ve) = 1 - \mathds{P}_{\tau_n,X_{\tau_n}}(S_{\sigma_n}=\ell^2+\ve)
\leq \mathds{P}_{\tau_n,X_{\tau_n}}(|G_n| \geq \ve/\sqrt{v_n}).
\end{equation*}
By the usual estimate of tails for the standard Gaussian random variables, we get (\ref{maj}).

{\sl  Step d).} To insure the convergence of the series in (\ref{but2}) we show that the sequence $(\tau_n)$ 
increases to infinity sufficiently fast. More precisely, we show that there exists $\lambda>1$ such that
$\tau_n\geq\lambda^n\,\tau_1 $ a.s. on $F$.
This inequality will be a consequence of a sharper form of the Borel-Cantelli lemma (see, for 
example, Theorem 1 in \cite{DF}, p. 800) once we show that there exist some constants $q>1$ and $p>0$ such that for all  $n\geq 1$, 
\begin{equation}\label{temps2}
\mathds{P}(\tau_{n+1}\geq q\,\tau_{n}\mid \mathcal{F}_{\tau_n})\geq \mathds{P}_{\tau_n,X_{\tau_n}}(\sigma_n\geq q\,\tau_n) 
\geq p
\quad\mbox{a.s}\quad\mbox{on}\quad\{\tau_n<\infty\}.
\end{equation}
In opposite to (\ref{estim_LV}), we can see that there exists a constant $k>0$ such that for all 
$t\geq 1$ and $x\in\mathds{R}$, for which  $G(t,x)\leq \ell^2+3\ve$, 
$LG(t,x)\geq {-k}/{t}$. We deduce that for all $t\in[\tau_n,q\tau_{n}]$,
\begin{equation*}
-\frac{\ve}{2}\leq k\ln{\left(\frac{\tau_n}{t}\right)}\leq A^n_t\leq 0,\quad \mbox{with}\quad q:=e^{\frac{k\ve }{2}}>1.
\end{equation*}
By using this inequality, we can write
\begin{equation*}
\left\{\sigma_n \geq  q \tau_n\right\}
\supset
\Big\{\sup_{\tau_n\leq t\leq q\tau_n} |M^n_t|<\frac{\varepsilon}{2}\Big\}
\supset
\Big\{\sup_{0\leq t\leq  v_0} |W^n_t|<\frac{\varepsilon}
{2}\Big\}\,.
\end{equation*}
Therefore, inequality (\ref{temps2}) is satisfied with the deterministic positive constant 
\begin{equation*}
p:=\mathds{P}_{\tau_n,X_{\tau_n}}\Big(\sup_{0\leq t\leq  v_0} |W^n_t|<\frac{\varepsilon}
{2}\Big)=\mathds{P}\Big(\sup_{0\leq t\leq  v_0} |B_t|<\frac{\varepsilon}{2}\Big).
\end{equation*}
Here $B$ denotes a standard Brownian motion. The sharper form of the Borel-Cantelli lemma applies 
and we obtain that $\tau_n\geq\lambda^n\,\tau_1 $ a.s. on $F$. We deduce (\ref{but2}) and  
(\ref{transient_scale}) holds.

Assume that $\alpha\in(1,\infty)$. The proof of (\ref{vitesseexplosion2}) follows the same lines as the proof 
of (\ref{vitesseexplosion}) in Theorem \ref{criticallinerepuls}. We show that $X^{(\upgamma)}=\Phi_\upgamma(X)$ satisfies
\begin{equation*}
|X_s^{(\upgamma)}|\underset{s \to\infty}{\sim}\frac{1}{(\rho(\alpha-1)(\tau_e-s) )^{\frac{1}{\alpha-1}}}\quad\mbox{a.s.,}\quad\mbox{when}\quad \alpha\in(1,\infty),
\end{equation*}
and we conclude by applying the change of time $t=\varphi^{-1}_\upgamma(s)$.

Finally, assume that $\alpha=1$. The same ideas as for the proof of the point v) of Theorem  (\ref{criticallinerepuls}) are employed. 
By Ito's formula and the Dambis-Dubins-Schwartz theorem there exists a standard Brownian motion $W$ such that
\begin{equation*}
\frac{X_t}{v(t)}= \frac{x_0}{v(1)}+{\tilde B}_{\phi(t)},\quad\mbox{with}\quad \phi(t):=\int_1^t \frac{ds}{v(s)^2}ds\quad\mbox{and}\quad v(t):=\exp\Big({\rho\frac{t^{1-\beta}}{1-\beta}}\Big).
\end{equation*}
Since
\begin{equation*}
 \phi(t)\underset{t\to\infty}{\sim}
\left\{
\begin{array}{cc}
t &\mbox{ if }\beta\in(1,\infty)\\
\sigma^2 &\mbox{ if }\rho\in(0,\infty)\mbox{ and }\beta\in(-\infty,1)\\
2|\rho|t^\beta v(t)^{-2} &\mbox{ if }\rho\in(-\infty,0)\mbox{ and }\beta\in(-\infty,1),
\end{array}
\right.
\end{equation*}
by using the usual properties of the Brownian motion we can get the convergence in distribution and the pathwise largest deviations. The  recurrent or transient features are then deduced. 
\hfill$\Box$

\section{Appendix}\label{lemfin}
\setcounter{equation}{0}

{\bf Proof of Lemma \ref{lemmepresquehorrible}.} To begin with, let us recall that $\rho\in(0,\infty)$, $\alpha\in(1,\infty)$, $2\beta\in(\alpha+1,\infty)$, that $\gamma=2\beta/(\alpha+1)$, $t_1=1/(\gamma-1)$, $\delta=\gamma/2(\gamma-1)$,
and that $\Hm_s$ is the pathwise unique strong solution of equation (\ref{eqtimechangegamma2}), which explosion 
time is $\eta_{e}\in[0,t_{1}]\cup\{\infty\}$. The goal is to prove that $\eta_{e}\neq t_{1}$ a.s. By Ito's formula, we can see that
\begin{equation*}
dX^{(\upgamma,t_1)}_s= (t_1-s)^{\frac{1}{\alpha-1}} dW_s + d(s,X^{(\upgamma,t_1)}_s)ds,\quad \mbox{with}\quad X^{(\upgamma,t_1)}_s:=(t_1-s)^{\frac{1}{\alpha-1}}\Hm_s,
\end{equation*}
and
\begin{equation*}
d(s,x):=\rho\frac{x\left(|x|^{\alpha-1}-\ell^{\alpha-1}\right)}{t_1-s}
\quad\mbox{and}\quad
\ell:=\Big(\frac{1+\delta(\alpha-1)}{\rho(\alpha-1)}\Big)^{\frac{1}{\alpha-1}}\in(0,\infty).
\end{equation*}
Roughly speaking, since $x\cdot d(s,x)\geq 0$ (respectively $\leq0$), according as $|x|\geq \ell$ (respectively $\leq \ell$), 
$0,-\infty$ and $\infty$ are ``attractive'' levels, whereas $-\ell$ and $\ell$ are ``repulsive'' levels for the process $X^{(\upgamma,t_1)}$.
The strategy of the proof is as follows: firstly, we show that 
\begin{equation}\label{lemm1}
\lim_{s\to t_1}\,|X^{(\upgamma,t_1)}_s|\in\{0,\ell,\infty\}\quad\mbox{a.s. on}\quad F:=\{\eta_e=t_1\}.
\end{equation}  
Secondly, we shall prove that the following three events are of probability zero:
\begin{multline}\label{trio}
F_{0}:=F\cap\Big\{\lim_{s\to t_1}\,|X^{(\upgamma,t_1)}_s|=0\Big\},\quad
F_{\ell}:=F\cap\Big\{\lim_{s\to t_1}\,|X^{(\upgamma,t_1)}_s|=\ell\Big\}\quad\mbox{ and,}\quad\\
F_{\infty}:=F\cap\Big\{\lim_{s\to t_1}\,|X^{(\upgamma,t_1)}_s|=\infty\Big\}.
\end{multline} 
We stress that the following reasoning will be performed by taking place on the event $F$. For simplicity, 
this condition will be understood and will dropped along the following five steps.

{\sl Step a).} We verify (\ref{lemm1}). Introduce 
$E:=\big\{\omega\in F:\liminf_{s\to t_1}X^{(\upgamma,t_1)}_s<\limsup_{s\to t_1}X^{(\upgamma,t_1)}_s\big\}$. 
Fix $\omega\in E$ and suppose that $\limsup_{s\to t_1}X^{(\upgamma,t_1)}_s(\omega)>\ell$. 
We can pick two sequences of real numbers (which depends on $\omega$), $(s_n)$ and $(u_n)$, 
such that $0\leq u_n\leq s_n<t_1$ for all integers $n$, $\lim_{n\to\infty}u_n=t_1$, and 
\begin{equation*}
X^{(\upgamma,t_1)}_{u_n}(\omega)-X^{(\upgamma,t_1)}_{s_n}(\omega)=\frac{1}{2}\Big(\limsup_{s\to t_1}X^{(\upgamma,t_1)}_s (\omega)-\ell\Big)>0.
\end{equation*}
Moreover, this choice could be done such that for any $s\in[u_n,s_n]$, $X^{(\upgamma,t_1)}_s(\omega)\geq \ell$.
Denote by $M$ the martingale part of $X^{(\upgamma,t_1)}$. Since the drift $d(s,x)$ is nonnegative for all 
$x\geq \ell$, we deduce 
\begin{multline*}
\left|M_{s_n}(\omega)-M_{u_n}(\omega)\right|
=
\Big|X^{(\upgamma,t_1)}_{u_n}(\omega)-X^{(\upgamma,t_1)}_{s_n}(\omega)+\int_{u_n}^{s_n}d(s,X^{(\upgamma,t_1)}_s(\omega)) ds\Big|\\
\geq \frac{1}{2}\Big(\limsup_{s\to t_1}X^{(\upgamma,t_1)}_s(\omega) -\ell\Big)>0.
\end{multline*}
A similar argument works when $0<\limsup_{s\to t_1}X^{(\upgamma,t_1)}_s(\omega)\leq \ell$, 
but also for the two symmetric situations $\liminf_{s\to t_1}X^{(\upgamma,t_1)}_s(\omega)<-\ell$ and 
$-\ell\leq\limsup_{s\to t_1}X^{(\upgamma,t_1)}_s(\omega)<0$. This means that 
\begin{equation*}
 \liminf_{n\to\infty}\left|M_{s_n}(\omega)-M_{u_n}(\omega)\right|>0,\quad\mbox{a.s.}\quad\mbox{on}\quad E.
\end{equation*}
Since $M$ is a.s. uniformly continuous on $[0,t_1]$, necessarily $\mathds P(E)=0$. We obtain equality (\ref{lemm1}) 
by noting that
\begin{equation*}
\lim_{s\to t_1}|d(s,X^{(\upgamma,t_1)}_s)|=\infty\quad\mbox{a.s. on }\quad F\cap\Big\{\lim_{s\to t_1}|X^{(\upgamma,t_1)}_s|\notin\{0,\ell\}\Big\}.
\end{equation*}

{\sl Step b).} Note that $h_s:=\Hm_s-W_s$ is the solution of the ordinary differential equation
\begin{equation*}\label{lemm5}
h^\prime_s = \rho\,{\rm sgn}(h_s+W_s)|h_s+W_s|^{\alpha}-\delta\,\frac{h_s+W_s}{t_1-s}.
\end{equation*}
We re-write the latter equation 
\begin{equation}\label{lemm6}
h_s^\prime = -\epsilon_1(X^{(\upgamma,t_1)}_s)\cdot\frac{\delta(h_s+W_s)}{t_1-s},\quad\mbox{with}\quad 
\ve_1(x):= 1-\frac{\rho}{\delta} |x|^{\alpha-1},
\end{equation}
and
\begin{equation}\label{lemm7}
{h_s^\prime} = \epsilon_2(X^{(\upgamma,t_1)}_s)\cdot\rho\,{\rm sgn}(h_s+W_s)|h_s+W_s|^\alpha,\quad\mbox{with}\quad
\ve_2(x):= 1-\frac{\delta}{\rho}|x|^{1-\alpha}.
\end{equation}

{\sl Step c).}  Recall that $\eta_{e}$ is the explosion time of $\Hm$. If we prove that $\Hm$ is bounded on $[0,t_{1}]$, a.s. on $F_0$, necessarily  $\mathds P(F_{0})=0$. Since $W$ is a.s. continuous on the compact 
$[0,t_{1}]$, it suffices to prove that $h$ is bounded on 
$[0,t_{1}]$, a.s. on $F_0$. Set $\kappa:=\sup_{s\in[0,t_{1}]}|W_s|$. We note that 
$\lim_{s\to t_1} \ve_1(X^{(\upgamma,t_1)}_s)=1$ a.s. on $F_0$. Therefore, for any $\omega\in F_{0}$, there 
exists $u\in [0,t_1)$ such that, for all $s\in[u,t_{1})$, $h_s(\omega)h_s^\prime(\omega)\mathds{1}_{\{|h_s(\omega)|\geq\kappa\}}\leq 0$, by using (\ref{lemm6}).
This implies that $h^2_s(\omega)$ is bounded on $[0,t_{1}]$ and we are done.

{\sl Step d).} If we prove that $\lim_{s\to t_1} |X^{(\upgamma,t_1)}_s| =(\rho(\alpha-1))^{\frac{1}{1-\alpha}}$ a.s. 
on $F_\infty$, then, necessarily $\mathds P(F_\infty)=0$. Clearly, $\lim_{s\to t_1} \ve_2(X^{(\upgamma,t_1)}_s)=1$ a.s. on $F_\infty$. Then, by using (\ref{lemm7}) and the fact that $W$ is bounded on $[0,t_1]$,
\begin{equation*}
\frac{|h_s|^{1-\alpha}}{\alpha-1}=\int_{s}^{t_1}\frac{h^\prime_u}{{\rm sgn}(h_u)|h_u|^\alpha}du\underset{s\to t_1}{\sim} \rho(t_1-s),\quad\mbox{a.s. on }F_\infty.
\end{equation*}
To conclude, it suffices to recall that $X^{(\upgamma,t_1)}_s=(t_1-s)^{\frac{1}{\alpha-1}}(h_{s}+W_{s})$.

{\sl Step e).} Similarly, if we prove that $\lim_{s\to t_1} |X^{(\upgamma,t_1)}_s| =\infty$ a.s. on $F_\ell$, then, necessarily $\mathds P(F_\ell)=0$. First, we show that 
\begin{equation}\label{lemm8}
 \lim_{s\to t_1}(X^{(\upgamma,t_1)}_s-\ell)^{2}=\infty,
\quad\mbox{ a.s. on }F_\ell^{+}:=F_\ell\cap\Big\{\lim_{s\to t_1}X^{(\upgamma,t_1)}_s=\ell\Big\}.
\end{equation}
Denote $K_s:=(X^{(\upgamma,t_1)}_s-\ell)^2$. By Ito's formula, we can write
\begin{equation*}
dK_s = 2(t_1-s)^{\frac{1}{\alpha-1}}\sqrt{K_s}\,dB_s +\Big(2q(X^{(\upgamma,t_1)}_s)\frac{K_s}{t_1-s}+\frac{1}{2}(t_1-s)^{\frac{2}{\alpha-1}} \Big) ds,
\end{equation*} 
where
\begin{equation*}
B_s:=\int_0^s{\rm sgn}(K_u)dW_u\quad\mbox{ and }\quad q(x):= \rho \,x\frac{|x|^{\alpha-1}-\ell^{\alpha-1}}{x-l}.
\end{equation*}
Introduce, for $v\in[0,t_1)$, $C_s(v)$ the pathwise unique strong solution of 
\begin{equation*}
dC_s(v) = 2(t_1-s)^{\frac{1}{\alpha-1}} \sqrt{C_s(v) }\,d{B}_s+\Big(q_\infty\frac{C_s(v)}{t_1-s}+\frac{1}{2}(t_1-s)^{\frac{2}{\alpha-1}}\Big)ds,\quad C_v(v) =K_v \mathds{1}_{\{\eta_e>v\}},
\end{equation*}
where $q_\infty:=\lim_{x\to\ell}q(x)=\rho(\alpha-1)\ell^{\alpha-1}$. 
By comparison and localisation (see also the proof of (\ref{comparison1})), we can show that 
for all $s\in[v,t_1)$,
\begin{equation*}
K_s\geq C_s(v)\quad\mbox{a.s. on }\quad \Omega_v:=F_\ell^{+}
\cap\Big\{\inf_{v\in[0,t_1)}|2q(X^{(\upgamma,t_1)}_s)|\geq q_\infty\Big\}.
\end{equation*}
By Ito's formula, the law of the process $C(v)$ equals to the law of the square of the unique weak solution $Q(v)$ 
of the equation
\begin{equation*}
dQ_s(v)=(t_1-s)^{\frac{1}{\alpha-1}}dB_s + \frac{q_\infty}{2}\frac{Q_s(v) }{t_1-s}\,ds,\quad Q_v(v)=\sqrt{C_v(v)}.
\end{equation*}
Since $Q(v)$ is the solution of a linear equation, it is not difficult to see that $\lim_{s\to t_1}|Q_s(v)|=\infty$ a.s. and then we deduce that $\lim_{s\to t_1}C_s(v)=\infty$ a.s. 
Hence, for any $v\in[0,t_1)$, $\lim_{s\to t_1}K_s=\infty$ a.s. on $\Omega_v$. Since $
\lim_{s\to t_1}2q(X^{(\upgamma,t_1)}_s)>q_\infty$ a.s. on $F_\ell^{+}$ we obtain that 
$\lim_{s\to t_1}K_s=\infty$ a.s. on $\cup_{v\in[0,t_1)}\Omega_v=F_\ell^{+}$, which is (\ref{lemm8}). 
We conclude that $\mathds P(F_\ell\cap\{\lim_{s\to t_1}X^{(\upgamma,t_1)}_s=\ell\})=0$.
Clearly by similar arguments, we can prove that $\mathds P(F_\ell\cap\{\lim_{s\to t_1}X^{(\upgamma,t_1)}_s=-\ell\})=0$.
Hence $\mathds P(F_\ell)=0$. The proof of the lemma is now complete.\hfill$\Box$\\

\noindent
{\bf Proof of Lemma \ref{convergence_distribution}.}
Denote by $\mathds{P}_{u,z}$ the distribution of the diffusion
$Z$ with $Z_u=z$ and
$\{{\rm T}_{u,s} :  0\leq u\leq s\}$ the associated time-inhomogeneous semi-group. Similarly, denote by $\mathds{P}_{z}$ the distribution of the diffusion $H$ starting from $z$ at initial time and $\left\{{\rm T}_s : s\geq 0 \right\}$ the associated semi-group. Clearly, the diffusion coefficient $(s,z)\mapsto a(u+s,z)$ and the drift $(s,z)\mapsto d(u+s,z)$ of the 
diffusion $s\mapsto Z_{u+s}$ satisfy the hypothesis of Theorem 11.1.4 in \cite{SV}, p. 264. We deduce 
that, for every $f\in{\rm C}_b([0,\infty);\mathds{R})$ and $s\in[ 0,\infty)$,
\begin{equation}\label{asymptergo3}
\lim_{u\to\infty}{\rm T}_{u,u+s}f(z)={\rm T}_sf(z)\quad
\mbox{ uniformly in }z\mbox{ on compact subsets of }\mathds{R}.
\end{equation}
Moreover,
\begin{equation}\label{asymptergo2}
\lim_{s\to\infty}{\rm T}_sf(z)=\Pi(f)
\quad \mbox{ uniformly in }z \mbox{ on compact subsets of } \mathds{R}.
\end{equation}
Indeed, assume that $z$ belongs to the compact set $[b,c]$. By using the strong Markov property, we can prove  that for all $s\in[ 0,\infty)$ and $v\in\mathds{R}$,
\begin{equation*}
\mathds{P}_{b}(H_s>v)\leq\mathds{P}_{z}(H_s>v)\leq\mathds{P}_{c}(H_s>v).
\end{equation*}
By using the ergodic theorem and these last inequalities we get the uniform convergence on compact subsets of $\mathds{R}$ in (\ref{asymptergo2}). 
Besides, by the Markov property, for all $s,u\in[ 0,\infty)$, 
\begin{equation*}
{\rm T}_{0,u+s}f(z_0) -\Pi(f) 
= {\rm T}_{0,u}\left[{\rm T}_{u,u+s}f- {\rm T}_sf \right] (z_0)
+ {\rm T}_{0,u}\left[{\rm T}_sf-\Pi(f) \right](z_0)
\end{equation*}
and clearly, for arbitrary $r,s,u$ nonnegative real numbers, 
\begin{equation*}
|{\rm T}_{0,u+s}f(z_0) -\Pi(f)|
\leq\sup_{z\in[-r,r]}\{|{\rm T}_{u,u+s}f(z)-{\rm T}_sf(z)|+
|{\rm T}_sf(z)-\Pi(f)|\}+4\|f\|_{\infty}\mathds{P}_{z_0}(|Z_u|\geq r).
\end{equation*}
Thanks to (\ref{asymptergo3}) and (\ref{asymptergo2}), for all $r,\ve>0$ there exists $s_0,u_0\in[ 0,\infty)$ such that for all $u\geq u_0$, 
\begin{equation*}
|{\rm T}_{0,u+s_0}f(z_0) -\Pi(f)|\leq \varepsilon + 4\|f\|_\infty \sup_{s\geq 0}\mathds{P}(|Z_s|\geq r).
\end{equation*}
Since $Z$ is bounded in probability we deduce that $\lim_{u\to\infty}{\rm T}_{0,u}f(z_{0})=\Pi(f)$.
\hfill$\Box$\\

\noindent
 {\bf Acknowledgments}
       
 The authors are grateful to the Referee and the Associate Editor for careful reading
 of the first version of the manuscript and for useful comments and remarks.

\end{document}